\documentclass[10pt]{article} % For LaTeX2e
\usepackage{fullpage}
\usepackage{amsmath,amssymb,amsfonts,amsthm}
\usepackage{graphicx}
\usepackage{color}
\usepackage{natbib}
\usepackage{verbatim}
\usepackage{url}
\usepackage[title]{appendix}

\title{Coordinate Descent Converges Faster with the\\ Gauss-Southwell Rule Than Random Selection}
\author{
Julie Nutini$^1$, Mark Schmidt$^1$, Issam H. Laradji$^1$, Michael Friedlander$^2$, Hoyt Koepke$^3$\\\\
$^1$University of British Columbia, $^2$University of California, Davis, $^3$Dato}
\date{}

\def\norm#1{\|#1\|}

\newtheorem{problem}{Problem}

\newcommand{\R}{{\rm I\!R}}

\newcommand{\diag}[1]{\hbox{diag}({#1})}
\newcommand{\prox}{\mathop{\rm prox}}
\newcommand{\argmin}[1]{\mathop{\hbox{argmin}}_{#1}}
\newcommand{\argmax}[1]{\mathop{\hbox{argmax}}_{#1}}
\newcommand{\sign}{\mathop{\rm sign}}

\newcommand{\sgn}[1]{\hbox{sgn}(#1)}

\definecolor{red}{rgb}{1,0,0}

\begin{document}
\maketitle

%%%%%%%%%%%%%%%%%%%%%%%%%%%%%%%%%%%%%%%%%%%%%%%%%%%

\begin{abstract} 
There has been significant recent work on the theory and application of randomized coordinate descent algorithms, beginning with the work of~\citeauthor{nesterov2012efficiency} [\textit{SIAM J. Optim., 22(2), 2012}], who showed that a random-coordinate selection rule achieves the same convergence rate as the Gauss-Southwell selection rule. This result
suggests that we should never use the Gauss-Southwell rule, because it is typically much more expensive than random selection. However, the empirical behaviours of these algorithms contradict this theoretical result: in applications where the computational costs of the selection rules are comparable, the Gauss-Southwell selection rule tends to perform substantially better than random coordinate selection. We give a simple analysis of the Gauss-Southwell rule showing that---except in extreme cases---its convergence rate is faster than choosing random coordinates. We also (i) show that exact coordinate optimization improves the convergence rate for certain sparse problems, (ii) propose a Gauss-Southwell-Lipschitz 
rule that gives an even faster convergence rate given knowledge of the Lipschitz constants of the partial derivatives, (iii) analyze the effect of approximate Gauss-Southwell rules, and (iv) analyze proximal-gradient variants of the Gauss-Southwell rule.
\end{abstract}

\section{Coordinate Descent Methods}

There has been substantial recent interest in applying coordinate
descent methods to solve large-scale optimization problems, starting with the seminal work
of~\citet{nesterov2012efficiency}, who gave the
first global rate-of-convergence analysis for coordinate-descent
methods for minimizing convex functions. This analysis suggests that
choosing a random coordinate to update gives the same performance
as choosing the ``best'' coordinate to update via the more expensive
Gauss-Southwell (GS) rule. (Nesterov also proposed a more clever
randomized scheme, which we consider later in this paper.) This result
gives a compelling argument to use randomized coordinate descent in
contexts where the GS rule is too expensive.
%, which has lead to large performance improvements on a variety of problems. 
It also suggests that there is no benefit to using the GS rule 
in contexts where it is relatively cheap. But in these contexts, the GS rule often 
substantially outperforms randomized coordinate selection in practice. 
This suggests that either the analysis of GS is not tight, or that
 there exists a class of functions for which the GS rule is as slow as 
randomized coordinate descent.

After discussing contexts in which it 
makes sense to use coordinate descent and the GS rule, we answer this
theoretical question by giving a tighter analysis of the GS rule (under 
strong-convexity and standard smoothness assumptions) that yields 
the same rate as the randomized method for a restricted class of functions, 
but is otherwise faster (and in some cases substantially faster).
We further show that, compared to the usual \emph{constant} step-size update 
of the coordinate, the GS method with exact coordinate optimization has 
a provably faster rate for problems satisfying a certain sparsity constraint 
(Section~\ref{sec:exact}). We believe that this is the first result showing 
a theoretical benefit of exact coordinate optimization; all previous analyses 
show that these strategies obtain the same rate as constant step-size updates, 
even though exact optimization tends to be faster in practice. Furthermore, in 
Section~\ref{sec:GSL}, we propose a variant of the GS rule 
that, similar to Nesterov's more clever randomized sampling scheme, uses 
knowledge of the Lipschitz constants of the coordinate-wise gradients to 
obtain a faster rate. We also analyze approximate GS rules 
(Section~\ref{sec:approx}), which provide an intermediate strategy
between randomized methods and the exact GS rule. Finally, we 
analyze proximal-gradient variants of the GS rule (Section~\ref{sec:prox}) for 
optimizing problems that include a separable non-smooth term.

\section{Problems of Interest}
\label{sec:gs-problems}

The rates of Nesterov show that coordinate descent can be faster than gradient descent in cases where, if we are optimizing $n$ variables,
the cost of performing $n$ coordinate updates is similar to the cost of 
performing one full gradient iteration. This essentially means that coordinate 
descent methods are useful for minimizing convex functions that can 
be expressed in one of the following two forms:
\begin{align*}
h_1(x)  := \sum_{i=1}^n g_i(x_i) + f(Ax),\qquad
h_2(x)  := \sum_{i \in V}g_i(x_i) + \sum_{(i,j)\in E}f_{ij}(x_i,x_j),
\end{align*}
where $x_i$ is element $i$ of $x$, $f$ is smooth and cheap, the $f_{ij}$ are smooth, $G=\{V,E\}$ is a
graph, and $A$ is a matrix. (It is assumed that all functions are
convex.)\footnote{We could also consider slightly more general cases like functions that are defined on hyper-edges~\citep{richtarik2012parallel}, provided that
we can still perform $n$ coordinate updates 
for a similar cost to one gradient evaluation.} The family of functions $h_1$
includes core machine-learning problems such as least squares, logistic
regression, lasso, and SVMs (when solved in dual
form)~\citep{hsieh2008dual}. Family $h_2$ includes quadratic functions, graph-based label propagation algorithms for
semi-supervised learning~\citep{bengio2006label}, and finding the most 
likely assignments in continuous pairwise graphical models~\citep{rue2005gaussian}.

% New version
In general, the GS rule for problem $h_2$ is as expensive as a full gradient evaluation. However, the structure of $G$ often allows efficient implementation of the GS rule. For example, if each node has at most $d$ neighbours, we can track the gradients of all 
the variables and use a max-heap structure to implement the GS rule in $O(d\log n)$ time~\citep{meshi2012convergence}.
This is similar to the cost of the randomized algorithm if $d \approx |E|/n$ (since the average cost of the randomized method depends on the average degree).
This condition is true in a variety of applications. For example, in spatial statistics we often use two-dimensional grid-structured graphs, where the maximum degree is four
and the average degree is slightly less than $4$.
As another example, for applying graph-based label 
propagation on the Facebook graph (to detect the spread of diseases, for example), the average number of friends is around $200$
but no user has more than seven thousand friends.\footnote{\url{https://recordsetter.com/world-record/facebook-friends}}
The maximum number of friends would be even smaller if we removed edges based on proximity.
A non-sparse example where GS is efficient is complete graphs, since here the average degree and maximum degree are both $(n-1)$.
Thus, the GS rule is efficient for optimizing dense quadratic functions. On the other hand, GS could be very inefficient for star graphs.

If each column of $A$ has at most $c$ non-zeroes and each row has at most $r$ non-zeroes, then for many notable instances of problem $h_1$ we can implement the GS rule in $O(cr\log n)$ time by maintaining $Ax$ as well as the gradient and again using a max-heap (see Appendix~\ref{app:h12}). Thus, GS will be efficient if $cr$ is similar to 
the number of non-zeroes in $A$ divided by $n$. Otherwise,~\citet{dhillon2011nearest} show that we can 
approximate the GS rule for problem $h_1$ with no $g_i$ functions by 
solving a nearest-neighbour problem. Their analysis of the GS rule in the convex case, 
however, gives the same convergence rate that is obtained by random selection (although the constant factor can be smaller by a factor of up to $n$). 
More recently,~\citet{shrivastava2014asymmetric} give a general method for 
approximating the GS rule for problem $h_1$ with no $g_i$ functions by writing it as a maximum 
inner-product search problem.

\section{Existing Analysis}

We are interested in solving the convex optimization problem
\begin{equation}\label{eq:opt-problem}
\min_{x \in \mathbb{R}^n}\ f(x),
\end{equation}
where $\nabla f$ is coordinate-wise $L$-Lipschitz continuous, i.e.,  
for each $i=1,\ldots,n$,
\[
|\nabla_i f(x+\alpha e_i) - \nabla_i f(x)| \leq L|\alpha|,
\quad\mbox{$\forall x\in\mathbb{R}^n$ and $\alpha\in\mathbb{R}$},
\]
where $e_i$ is a vector with a one in position $i$ and zero in all other positions.
%where $e_i$ is a vector of zeros with a $1$ in the $i$th coordinate.
For twice-differentiable functions, this is equivalent to the
assumption that the diagonal elements of the Hessian are bounded in
magnitude by $L$. In contrast, the typical assumption used for
gradient methods is that $\nabla f$ is $L^f$-Lipschitz continuous (note that $L \leq L^f \leq Ln$). The coordinate-descent method with constant
step-size is based on the iteration
\[
x^{k+1} = x^k - \frac{1}{L}\nabla_{i_k} f(x^k)e_{i_k}.
\]
The randomized coordinate-selection rule chooses $i_k$ uniformly from
the set $\{1,2,\dots,n\}$. Alternatively, the GS rule
\[
i_k = \argmax{i}\ |\nabla_i f(x^k)|,
\]
chooses the coordinate with the largest directional derivative. Under
either rule, because $f$ is coordinate-wise Lipschitz continuous, we
obtain the following bound on the progress made by each iteration:
\begin{equation} \label{eq:descent}
\begin{aligned}
 f(x^{k+1})
& \leq f(x^k) + \nabla_{i_k}f(x^k)(x^{k+1}-x^k)_{i_k} + \frac{L}{2}(x^{k+1} - x^k)_{i_k}^2\\
& = f(x^k) -\frac1L (\nabla_{i_k}f(x^k))^2 + \frac{L}{2}\left[\frac{1}{L}\nabla_{i_k}f(x^k)\right]^2 \\
& = f(x^k) - \frac{1}{2L}[\nabla_{i_k}f(x^k)]^2.
\end{aligned}
\end{equation}
We focus on the case where $f$ is $\mu$-strongly convex, meaning that, 
for some positive $\mu$,
\begin{equation}\label{eq:strong-convexity1}
f(y) \geq f(x) + \langle \nabla f(x),y-x\rangle + \frac{\mu}{2}\norm{y-x}^2, \quad\forall x,y\in\mathbb{R}^n,
\end{equation}
which implies that
\begin{equation} \label{eq:strong-convexity}
f(x^*) \geq f(x^k) - \frac{1}{2\mu}\norm{\nabla f(x^k)}^2,
\end{equation}
where $x^*$ is the optimal solution of~\eqref{eq:opt-problem}. This
bound is obtained by minimizing both sides of~\eqref{eq:strong-convexity1} with respect to $y$.

\subsection{Randomized Coordinate Descent}

Conditioning on the $\sigma$-field $\mathcal{F}_{k-1}$ generated by the sequence
$\{x^0,x^1,\ldots,x^{k-1}\}$, and taking expectations of both sides of~\eqref{eq:descent},
when $i_k$ is chosen with uniform sampling we obtain
\begin{align*}
\mathbb{E}[f(x^{k+1})] 
& \leq \mathbb{E}\left[f(x^k) - \frac{1}{2L}\big(\nabla_{i_k}f(x^k)\big)^2\right]\\
& = f(x^k) - \frac{1}{2L}\sum_{i=1}^n \frac{1}{n}\big(\nabla_i f(x^k)\big)^2\\
& = f(x^k) - \frac{1}{2Ln}\norm{\nabla f(x^k)}^2.
\end{align*}
Using~\eqref{eq:strong-convexity} and subtracting $f(x^*)$ from both sides, we get
\begin{equation} \label{eq:random-descent}
\mathbb{E}[f(x^{k+1})] - f(x^*) \leq \left(1 - \frac{\mu}{Ln}\right)[f(x^k) - f(x^*)].
\end{equation}
This is a special of case of~\citet[][Theorem 2]{nesterov2012efficiency} with $\alpha=0$ in his notation.

\subsection{Gauss-Southwell}\label{sec:gauss-southwell}

We now consider the progress implied by the GS rule. By the definition of $i_k$,
\begin{equation} \label{eq:2-inf-inequality}
(\nabla_{i_k}f(x^k))^2 = \norm{\nabla f(x^k)}_\infty^2 \geq (1/n)\norm{\nabla f(x^k)}^2.
\end{equation}
Applying this inequality to~\eqref{eq:descent}, we obtain
\[
f(x^{k+1}) \leq f(x^k) - \frac{1}{2Ln}\norm{\nabla f(x^k)}^2,
\]
which together with~\eqref{eq:strong-convexity}, implies that
\begin{equation} \label{eq:gs-old}
f(x^{k+1}) - f(x^*) \leq \left(1 - \frac{\mu}{Ln}\right)[f(x^k) - f(x^*)].
\end{equation}
This is a special case of~\citet[][\S9.4.3]{boyd2004convex}, viewing the GS rule as performing steepest descent in the $1$-norm.
While this is faster than known rates for cyclic coordinate selection~\citep{beck2013convergence} and holds deterministically rather than in expectation, this rate is the same as the randomized rate given in~\eqref{eq:random-descent}.

\section{Refined Gauss-Southwell Analysis}
\label{sec:new}

The deficiency of the existing GS analysis is that too much is
lost when we use the inequality in~\eqref{eq:2-inf-inequality}. To
avoid the need to use this inequality, we instead measure strong-convexity in
the $1$-norm, i.e.,
\[
f(y) \geq f(x) + \langle \nabla f(x),y-x\rangle + \frac{\mu_1}{2}\norm{y-x}_1^2,
\]
which is the analogue of~\eqref{eq:strong-convexity1}. Minimizing both 
sides with respect to $y$, we obtain
\begin{equation}
\label{eq:LB}
\begin{aligned}
f(x^*) & \geq f(x) - \sup_y\{\langle -\nabla f(x),y-x\rangle - \frac{\mu_1}{2}\norm{y-x}_1^2\}\\
& = f(x) - \left(\frac{\mu_1}{2}\norm{\cdot}_1^2\right)^*(-\nabla f(x)) \\
& = f(x) - \frac{1}{2\mu_1}\norm{\nabla f(x)}_\infty^2,
\end{aligned}
\end{equation}
% If we minimize
% both sides of this inequality, we obtain
% \begin{align*}
% f(x^*) & \geq \inf_y\ \left\{ f(x) + \langle \nabla f(x),y - x\rangle  + \frac{\mu_1}{2}\norm{y-x}_1^2\right\}\\
% & = f(x) - \sup_y\ \{\langle -\nabla f(x),y-x\rangle - \frac{\mu_1}{2}\norm{y-x}_1^2\}\\
% & = f(x) - \left( \frac{\mu_1}{2}\norm{\cdot}_1^2 \right)^*\big({-}\nabla f(x)\big)\\
% & = f(x) - \frac{1}{2\mu_1}\norm{\nabla f(x)}_\infty^2,
% \end{align*}
% where we use the fact that the conjugate
% $\left(\frac{\mu_1}{2}\norm{\cdot}_1^2\right)^*=\frac1{2\mu_1}\norm{\cdot}_\infty^2$.
%Combine the previous two inequalities, with $x=x^k$ and $y=x^*$, to
%obtain
which makes use of the convex conjugate
$(\frac{\mu_1}{2}\norm{\cdot}_1^2)^* =
\frac{1}{2\mu_1}\norm{\cdot}_\infty^2$
\citep[\S3.3]{boyd2004convex}.  Using~\eqref{eq:LB}
in~\eqref{eq:descent}, and the fact that
$(\nabla_{i_k}f(x^k))^2 = \norm{\nabla f(x^k)}_\infty^2$ for the GS
rule, we obtain
\begin{equation} \label{eq:gs-new}
f(x^{k+1}) - f(x^*) \leq \left(1 - \frac{\mu_1}{L}\right)[f(x^k) - f(x^*)].
\end{equation}

%\subsection{Comparison of Random and Gauss-Southwell}
%\label{sec:mu12}

It is evident that if $\mu_1 = \mu/n$, then the rates implied
by~\eqref{eq:random-descent} and~\eqref{eq:gs-new} are identical, 
but~\eqref{eq:gs-new} is faster if $\mu_1 > \mu/n$.
%It follows from \citet[Theorem 2.1.9]{Nes04b} that the strong-convexity
%parameters $\mu$ and $\mu_1$ can be defined by
%\begin{align*}
%\mu & = \inf_{x,y}\ \frac{\norm{\nabla f(x)-\nabla f(y)}}{\norm{x-y}},\\
%\mu_1 &= \inf_{x,y}\ \frac{\norm{\nabla f(x)-\nabla f(y)}_\infty}{\norm{x-y}_1}.
%\end{align*}
In Appendix~\ref{app:mu12}, we show that the relationship between $\mu$ and $\mu_1$ can be obtained through the relationship between the squared norms $||\cdot||^2$ and $||\cdot||_1^2$. In particular, we have
\[
\frac{\mu}{n} \leq \mu_1 \leq \mu.
\]
Thus, at one extreme the GS rule obtains the same rate as uniform 
selection ($\mu_1 \approx \mu/n$). However, at the other extreme, it could be 
faster than uniform selection by a factor of $n$ ($\mu_1 \approx \mu$). This 
analysis, that the GS rule only obtains the same bound as random selection in 
an extreme case, supports  the better practical behaviour of GS.
%the GS rule having a better convergence rate estimate
%than random selection.
%why the GS rule often behaves much better than 
%random selection in practice.

\subsection{Comparison for Separable Quadratic}
\label{sec:mu1}

We illustrate these two extremes with the simple example of a
quadratic function with a diagonal Hessian $\nabla^2 f(x) =
\diag{\lambda_1,\ldots,\lambda_n}$. In this case,
\[
\mu = \min_i\ \lambda_i,
\quad\mbox{and}\quad
\mu_1 = \left(\sum_{i=1}^n\frac{1}{\lambda_i}\right)^{-1}.
%\frac{\prod_{i=1}^n \lambda_i}{\sum_{k=1}^n\prod_{i\ne k}\lambda_i}.
\]
We prove the correctness of this formula for $\mu_1$ in Appendix~\ref{app:mu1}.
The parameter $\mu_1$ achieves its lower bound when all $\lambda_i$ are equal,
$\lambda_1=\cdots=\lambda_n=\alpha>0$, in which case
\[
\mu = \alpha \quad\mbox{and}\quad \mu_1 = \alpha/n.
\]
Thus, uniform selection does as well as the GS rule if all elements of the 
gradient change at \emph{exactly} the same rate. This is reasonable: under this 
condition, there is no apparent advantage in selecting the coordinate to update 
in a clever way. Intuitively, one might expect that the favourable case for the 
Gauss-Southwell rule would be where one $\lambda_i$ is much larger than the 
others. However, in this case, $\mu_1$ is again similar to $\mu/n$.  To achieve the 
other extreme, suppose that $\lambda_1 = \beta$ and 
$\lambda_2=\lambda_3=\cdots=\lambda_n=\alpha$ with $\alpha \geq \beta$. 
In this case, we have $\mu = \beta$ and
\[
\mu_1 = \frac{\beta\alpha^{n-1}}{\alpha^{n-1}+ (n-1)\beta\alpha^{n-2}} = \frac{\beta\alpha}{\alpha + (n-1)\beta}.
\]
If we take $\alpha\to\infty$, then we have $\mu_1\to\beta$, so $\mu_1\to\mu$. 
This case is much less intuitive; GS is $n$ times faster than random coordinate 
selection if one element of the gradient changes much more \emph{slowly} 
than the others.
% Appendix~4.1
%gives a physical interpretation of $\mu$ and $\mu_1$
%in terms of independent processes `working together'~\cite{ferger1931harmonic}.

\subsection{`Working Together' Interpretation}

In the separable quadratic case above, $\mu_1$ is given by the harmonic mean
of the eigenvalues of the Hessian divided by $n$. The harmonic mean is dominated by its smallest values, and this is why having one small value is a notable case. Furthermore, the harmonic mean divided by $n$ has an interpretation in terms of processes `working
together' \citep{ferger1931harmonic}.  If each $\lambda_i$
represents the time taken by each process to finish a task (e.g.,
large values of $\lambda_i$ correspond to slow workers), then $\mu$ is
the time needed by the fastest worker to complete the task, and
$\mu_1$ is the time needed to complete the task if all processes work
together (and have independent effects).
Using this interpretation,
the GS rule provides the most benefit over random selection when
\emph{working together is not efficient}, meaning that if the $n$ processes work together, then the task is not solved much faster than if the fastest worker performed the task alone.  This gives an interpretation of the
non-intuitive scenario where GS provides the most benefit: if all
workers have the same efficiency, then working together solves the
problem $n$ times faster. Similarly, if there is one slow worker
(large $\lambda_i$), then the problem is solved roughly $n$ times
faster by working together. On the other hand, if most workers are
slow (many large $\lambda_i$), then working together has little
benefit.

\subsection{Fast Convergence with Bias Term}

Consider the standard linear-prediction framework,
\[
\argmin{x,\beta} \sum_{i=1}^m f(a_i^Tx + \beta) + \frac{\lambda}{2}\norm{x}^2 + \frac{\sigma}{2}\beta^2,
\]
where we have included a bias variable $\beta$ (an example of problem $h_1$). Typically, the regularization 
parameter $\sigma$ of the bias variable is set to be much smaller than the 
regularization parameter $\lambda$ of the other covariates, to avoid biasing 
against a global shift in the predictor. Assuming that there is no hidden 
strong-convexity in the sum, this problem has the structure described in the previous section ($\mu_1 \approx \mu$) where GS has the 
most benefit over random selection.

\section{Rates with Different Lipschitz Constants}
\label{sec:exact}

Consider the more general scenario where we have a Lipschitz constant 
$L_i$ for the partial derivative of $f$ with respect to each coordinate $i$,
\[
|\nabla_i f(x+\alpha e_i) - \nabla_i f(x)| \leq L_i|\alpha|,
\quad\mbox{$\forall x\in\mathbb{R}^n$ and $\alpha\in\mathbb{R}$,}
\]
and we use a coordinate-dependent step-size at each iteration:
\begin{equation}
\label{eq:CD_i}
x^{k+1} = x^k - \frac{1}{L_{i_k}}\nabla_{i_k} f(x^k)e_{i_k}.
\end{equation}
By the logic of~\eqref{eq:descent}, in this setting we have
\begin{equation}
\label{eq:descent_i}
f(x^{k+1}) \leq f(x^k) - \frac{1}{2L_{i_k}}[\nabla_{i_k}f(x^k)]^2,
\end{equation}
and thus a convergence rate of
\begin{equation}
\label{eq:ratei}
f(x^k) - f(x^*) \leq \left[\prod_{j=1}^k\left(1 - \frac{\mu_1}{L_{i_j}}\right)\right][f(x^0) - f(x^*)].
\end{equation}
Noting that $L = \max_i \{ L_i \}$, we have
\begin{equation}
\label{eq:Li_rate}
 \prod_{j=1}^k\left(1 - \frac{\mu_1}{L_{i_j}}\right) \leq \left(1-\frac{\mu_1}{L}\right)^k.
\end{equation}
Thus, the convergence rate based on the $L_i$ will be faster, provided that
at least one iteration chooses an $i_k$ with $L_{i_k} < L$. In the
worst case, however,~\eqref{eq:Li_rate} holds with equality even if
the $L_i$ are distinct, as we might need to update a coordinate with
$L_i = L$ on every iteration. (For example, consider a separable
function where all but one coordinate is initialized at its optimal
value, and the remaining coordinate has $L_i = L$.)  In
Section~\ref{sec:GSL}, we discuss selection rules that incorporate the
$L_i$ to achieve faster rates whenever the $L_i$ are distinct, but first
we consider the effect of exact coordinate optimization on the choice of the $L_{i_k}$.
%But first, we show that using 
%exact coordinate optimization guarantees a faster convergence rate when using 
%the Gauss-Southwell rule for sparse versions of problems $h_1$ and $h_2$.

\subsection{Gauss-Southwell with Exact Optimization}

For problems involving functions of the form $h_1$ and $h_2$, we are often 
able to perform exact (or numerically very precise) coordinate optimization, even 
if the objective function is not quadratic (e.g., by using a line-search or a closed-form 
update). Note that~\eqref{eq:ratei} still holds when using exact coordinate 
optimization rather than using a step-size of $1/L_{i_k}$, as in this case we have
\begin{equation}
\begin{aligned}
f(x^{k+1}) & = \min_{\alpha}\{f(x^k + \alpha e_{i_k})\}\\
& \leq f\left(x^k  -  \frac{1}{L_{i_k}}\nabla_{i_i}f(x^k)e_{i_k}\right) \\
& \leq f(x^k) - \frac{1}{2L_{i_k}}[\nabla_{i_k}f(x^k)]^2,
\end{aligned}
\label{eq:LS}
\end{equation}
which is equivalent to~\eqref{eq:descent_i}. 
%In the worst case, exact coordinate 
%optimization cannot improve the per-iteration convergence rate because it coincides 
%with using a step-size of $1/L_{i_k}$ in the special case of quadratic functions. 
However, in practice using exact coordinate optimization leads to better performance. In this section, we show that using the GS rule results in a convergence 
rate that is indeed faster than~\eqref{eq:gs-new} for problems with distinct $L_i$ when the function is 
quadratic, or when the function is not quadratic but we perform exact coordinate optimization.
%Although we can't show a faster convergence rate in general for exact coordinate optimization (because the two updates coincide for quadratic functions), in this section we show that coordinate optimization can be much faster for certain sparse problems.

The key property we use is that, after we have performed exact coordinate optimization, 
we are guaranteed to have $\nabla_{i_k}f(x^{k+1}) = 0$. Because the GS rule chooses 
$i_{k+1} = \argmax{i}|\nabla_i f(x^{k+1})|$, we cannot have $i_{k+1} = i_k$, unless $x^{k+1}$ 
is the optimal solution. Hence, we never choose the same coordinate twice in a 
row, which guarantees that the inequality~\eqref{eq:Li_rate} is strict (with distinct $L_i$) 
and exact coordinate optimization is faster. We note that the improvement may be marginal, 
as we may simply alternate between the two largest $L_i$ values. However, consider 
minimizing $h_2$ when the graph is sparse; after updating 
$i_k$, we are guaranteed to have $\nabla_{i_k}f(x^{k+m}) = 0$ for all future iterations $(k+m)$ 
until we choose a variable $i_{k+m-1}$ that is a neighbour of node $i_k$ in the graph. 
Thus, if the two largest $L_i$ are not connected in the graph, GS cannot 
simply alternate between the two largest $L_i$. 

By using this property, in Appendix~\ref{app:exact} we show that the GS rule with exact coordinate optimization for problem $h_2$ under a chain-structured graph has a convergence rate of the form
\begin{align*}
 f(x^k) - f(x^*) \leq O\left(\max\{\rho_2^G,\rho_3^G\}^k\right)[f(x^0) - f(x^*)],
\end{align*}
where $\rho_2^G$ is the maximizer of $\sqrt{(1-\mu_1/L_i)(1-\mu_1/L_j)}$ among all consecutive nodes $i$ and $j$ in the chain, and $\rho_3^G$ is the maximizer of $\sqrt[3]{(1-\mu_1/L_i)(1-\mu_1/L_j)(1-\mu_1/L_k)}$ among consecutive nodes $i$, $j$, and $k$.
%triplets $\{i,j,k\}$, where we have an edge between $i$ and $j$ and an edge between $j$ and $k$. 
The implication 
of this result is that, if the large $L_i$ values are more than two edges from each other in the graph, then we obtain a much better convergence rate. 
We conjecture that for general graphs, we can obtain a bound that depends on the largest value of $\rho_2^G$ among all nodes $i$ and $j$ connected by a path of length $1$ or $2$.
Note that we can obtain similar results for problem $h_1$, by forming a graph that has an edge between nodes $i$ and $j$ whenever the corresponding variables are both jointly non-zero in at least one  row of $A$.

\section{Rules Depending on Lipschitz Constants}
\label{sec:GSL}

If the $L_i$ are known,~\citet{nesterov2012efficiency} showed that we can obtain 
a faster convergence rate by sampling proportional to the $L_i$. We review this 
result below and compare it to the GS rule, and then propose an improved GS 
rule for this  scenario. Although in this section we will assume that the $L_i$ are known, 
this assumption can be relaxed using a backtracking procedure~\citep[][\S6.1]{nesterov2012efficiency}.

\subsection{Lipschitz Sampling}

Taking the expectation of~\eqref{eq:descent_i} under the distribution 
$p_i = L_i/\sum_{j=1}^nL_j$ and proceeding as before, we obtain
\[
\mathbb{E}[f(x^{k+1})] - f(x^*) \leq \left(1-\frac{\mu}{n\bar{L}}\right)[f(x^k)-f(x^*)],
\]
where $\bar{L}=\frac{1}{n}\sum_{j=1}^nL_j$ is the average of the Lipschitz constants.
This was shown by~\citet{leventhal2010randomized} and is a special case of~\citet[][Theorem 2]{nesterov2012efficiency} with $\alpha=1$ 
in his notation. This rate is faster than~\eqref{eq:random-descent} for uniform 
sampling if any $L_i$ differ.

Under our analysis, this rate may or may not be faster than~\eqref{eq:gs-new} 
for the GS rule. On the one extreme, if $\mu_1=\mu/n$ and any $L_i$ differ, then 
this Lipschitz sampling scheme is faster than our rate for GS. Indeed, in the context 
of the problem from Section~\ref{sec:mu1}, we can make Lipschitz sampling faster 
than GS by a factor of nearly $n$ by making one $\lambda_i$ much larger than 
all the others (recall that our analysis shows no benefit to the GS rule over randomized 
selection when only one $\lambda_i$ is much larger than the others). At the other extreme, 
in our example from Section~\ref{sec:mu1} with many large $\alpha$ and one small $\beta$, 
the GS and Lipschitz sampling rates are the same when $n=2$, with a rate of 
$(1-\beta/(\alpha+\beta))$. However, the GS rate will be faster than the Lipschitz 
sampling rate for any $\alpha > \beta$ when $n > 2$, as the Lipschitz sampling rate 
is $(1-\beta/((n-1)\alpha + \beta))$, which is slower than the GS rate of $(1-\beta/(\alpha + (n-1)\beta))$. 

\subsection{Gauss-Southwell-Lipschitz Rule}

Since neither Lipschitz sampling nor GS dominates the other in
general, we are motivated to 
consider if faster rules are possible by combining the two approaches. Indeed, 
we obtain a faster rate by choosing the $i_k$ that minimizes~\eqref{eq:descent_i}, leading 
to the rule
\[
i_k = \argmax{i} \frac{|\nabla_i f(x^k)|}{\sqrt{L_i}},
\]
which we call the \emph{Gauss-Southwell-Lipschitz} (GSL) rule. Following a similar argument to 
Section~\ref{sec:new}, but using~\eqref{eq:descent_i} in place of~\eqref{eq:descent}, the GSL rule 
obtains a convergence rate of
\[
f(x^{k+1}) - f(x^*) \leq (1-\mu_L)[f(x^k) - f(x^*)],
\]
where $\mu_L$ is the strong-convexity constant with respect to the norm 
$\norm{x}_L = \sum_{i=1}^n\sqrt{L_i}|x_i|$. This is shown in Appendix~\ref{app:GSL}, and in Appendix~\ref{app:muL}
we show that
\[
\max\left\{\frac{\mu}{n\bar{L}},\frac{\mu_1}{L}\right\} \leq \mu_L \leq \frac{\mu_1}{\min_i\{L_i\}}.
\]
Thus, the GSL rule is always at least as fast as the fastest of the GS rule 
and Lipschitz sampling. Indeed, it can be more than a factor of $n$ faster than using 
Lipschitz sampling, while it can obtain a rate closer to the minimum $L_i$, instead of the 
maximum $L_i$ that the classic GS rule depends on.
%If all $L_i$ are equal, then $\mu_L = \mu_1/L$ and we obtain~\eqref{eq:gs-new}. But otherwise, this refined GS rule gives a faster rate that be as fast as $\mu_1$ divided by the smallest $L_i$ rather than the largest value $L$. This result shows that the GSL rule can be substantially faster than Lipscthiz sampling for many problems, and we suspect that it's possible to improve the $\sqrt{L\bar{[L}}$ factor to $\bar{L}$ to show that the GSL rule is never slower than Lipschitz sampling.

An interesting property of the GSL rule for quadratic functions is that it is the \emph{optimal} myopic coordinate 
update. That is, if we have an oracle that can choose the 
coordinate and the step-size that decreases $f$ by the largest amount, i.e.,
\begin{equation}
\label{eq:maxImprov}
f(x^{k+1}) = \argmin{i,\alpha}\{f(x^k + \alpha e_i)\},
\end{equation}
this is equivalent to using the GSL rule and the update in~\eqref{eq:CD_i}.
This follows because~\eqref{eq:descent_i} holds with equality in the quadratic case, 
and the choice $\alpha_k = 1/L_{i_k}$ yields the optimal step-size. Thus, although 
faster schemes could be possible with non-myopic strategies that cleverly choose 
the sequence of coordinates or step-sizes, if we can only perform one iteration, 
then the GSL rule cannot be improved.

For general $f$,~\eqref{eq:maxImprov} is known as the \emph{maximum improvement} (MI) rule.
This rule has been used in the context of boosting~\citep{ratsch2001convergence}, graphical models~\citep{della1997inducing,lee2006efficient,scheinberg2009sinco}, Gaussian processes~\citep{bo2012greedy}, and low-rank tensor approximations~\citep{li2015convergence}.
Using an argument similar to~\eqref{eq:LS}, our GSL rate also applies to the MI rule, improving existing bounds on this strategy.
 However, the GSL rule is much cheaper and does not require any special structure (recall that we can estimate $L_i$ as we go). 
 
 \subsection{Connection between GSL Rule and Normalized Nearest Neighbour Search}
 \label{sec:nns}
 
\citet{dhillon2011nearest} discuss an interesting connection between the GS rule and the nearest-neighbour-search (NNS) problem for objectives of the form
\begin{equation}
\label{eq:basicH1}
\min_{x \in \R^n} F(x) = f(Ax),
\end{equation}
This is a special case of $h_1$ with no $g_i$ functions, and its gradient has the special form
\[
\nabla F(x) = A^T r(x),
\]
where $r(x) = \nabla f(Ax)$. We use the symbol $r$ because it is the residual vector ($Ax-b$) in the special case of least squares. For this problem structure the GS rule has the form
\begin{align*}
i_k & = \argmax{i}|\nabla_i f(x^k)| \\
& = \argmax{i}|r(x^k)^Ta_i|,
\end{align*}
%Using $\mathcal{A} = \{a_1, a_2, \dots, a_m, -a_1, -a_2, \dots -a_n \}$ be an ordered set containing the columns $a_i$ and their negations $a_i$,
where $a_i$ denotes column $i$ of $A$ for $i = 1, \dots, n$. \citet{dhillon2011nearest} propose to approximate the above $\argmax{}$ by solving the following NNS problem
\[
i_k = \argmin{i\in [2n]} \norm{r(x^k) - a_i},
\]
where $i$ in the range $(n+1)$ through $2n$ refers to the negation $-(a_{i-n})$ of column $(i-n)$ and if the selected $i_k$ is greater than $n$ we return $(i-n)$. We can justify this approximation using the logic
\begin{align*}
i_k & = \argmin{i\in [2n]} \norm{r(x^k)-a_i} \\
& = \argmin{i\in [2n]} \frac{1}{2}\norm{r(x^k)-a_i}^2 \\
& = \argmin{i\in [2n]}  \underbrace{\frac{1}{2}\norm{r(x^k)}^2}_{\text{constant}} - r(x^k)^Ta_i + \frac{1}{2}\norm{a_i}^2\\
%& = \argmin{i\in [2n]}  - r(x)^Ta_i + \frac{1}{2}\norm{a_i}^2\\
& = \argmax{i\in [2n]}  r(x^k)^Ta_i - \frac{1}{2}\norm{a_i}^2\\
& = \argmax{i\in [n]}  |r(x^k)^Ta_i| - \frac{1}{2}\norm{a_i}^2\\
& = \argmax{i\in [n]} |\nabla_i f(x^k)| - \frac{1}{2}\norm{a_i}^2.
\end{align*}
Thus, the NNS computes an approximation to the GS rule that is biased towards coordinates where $\norm{a_i}$ is small. Note that this formulation is equivalent to the GS rule in the special case that $\norm{a_i}=1$ (or any other constant) for all $i$. \citet{shrivastava2014asymmetric} have more recently considered the case where $\norm{a_i} \leq 1$ and incorporate powers of $\norm{a_i}$ in the NNS to yield a better approximation. 
 
 A further interesting property of the GSL rule is that we can often formulate the \emph{exact} GSL rule as a \emph{normalized} NNS problem. In particular, for problem~\eqref{eq:basicH1} the Lipschitz constants will often have the form $L_i = \gamma\norm{a_i}^2$ for a some positive scalar $\gamma$. For example, least squares has $\gamma = 1$ and logistic regression has $\gamma = 0.25$. When the Lipschitz constants have this form, we can compute the exact GSL rule by solving a normalized NNS problem,
 \begin{equation}
\label{eq:GSL_NNS}
i_k = \argmin{i \in [2n]} \left|\left|r(x^k) - \frac{a_i}{\norm{a_i}}\right|\right|.
\end{equation}
The exactness of this formula follows because
\begin{align*}
i_k 
& = \argmin{i \in [2n]} \left|\left|r(x^k) - \frac{a_i}{\norm{a_i}}\right|\right|\\
& = \argmin{i\in [2n]} \frac{1}{2}\norm{r(x^k)-a_i/\norm{a_i}}^2 \\
& = \argmin{i\in [2n]}  \underbrace{\frac{1}{2}\norm{r(x^k)}^2}_{\text{constant}} - \frac{r(x^k)^Ta_i}{\norm{a_i}} + \underbrace{\frac{1}{2}\frac{\norm{a_i}^2}{\norm{a_i}^2}}_{\text{constant}}\\
%& = \argmin{i\in [2n]}  - r(x)^Ta_i + \frac{1}{2}\norm{a_i}^2\\
%& = \argmax{i\in [2n]}  r(x)^Ta_i - \frac{1}{2}\norm{a_i}^2\\
& = \argmax{i\in [n]}  \frac{|r(x^k)^Ta_i|}{\norm{a_i}}\\
& = \argmax{i\in [n]}  \frac{|r(x^k)^Ta_i|}{\sqrt{\gamma}\norm{a_i}}\\
& = \argmax{i\in [n]} \frac{|\nabla_i f(x^k)|}{\sqrt{L_i}}.
\end{align*}
Thus, the form of the Lipschitz constant conveniently removes the bias towards smaller values of $\norm{a_i}$ that gets introduced when we try to formulate the classic GS rule as a NNS problem. Interestingly, in this setting we \emph{do not need to know $\gamma$} to implement the GSL rule as a NNS problem.

\section{Approximate Gauss-Southwell}
\label{sec:approx}

In many applications, computing the exact GS rule is too
inefficient to be of any practical use. However, a computationally cheaper \emph{approximate} 
GS rule might be available. 
%For example, we can approximate the GS rule for 
%functions of the form $h_1$ using fast approximate maximum inner-product search 
%methods~\citep{dhillon2011nearest,shrivastava2014asymmetric}.
Approximate GS rules 
under multiplicative and additive errors were considered by~\citet{dhillon2011nearest} in the convex case, but in this setting the convergence rate is similar to the rate achieved by random selection.
In this section, we give rates depending on $\mu_1$ for approximate GS rules.

\subsection{Multiplicative Errors}

In the multiplicative error regime, the approximate GS rule chooses an $i_k$ satisfying
\[
|\nabla_{i_k} f(x^k)| \geq \norm{\nabla f(x^k)}_\infty(1 - \epsilon_k),
\]
for some $\epsilon_k \in[0,1)$. In this regime, our basic bound on the progress~\eqref{eq:descent} 
still holds, as it was defined for any $i_k$. We can incorporate this type of error 
into our lower bound~\eqref{eq:LB} to obtain
\begin{align*}
f(x^*) & \geq f(x^k) - \frac{1}{2\mu_1}\norm{\nabla f(x^k)}_\infty^2\\
& \geq f(x^k) - \frac{1}{2\mu_1(1-\epsilon_k)^2}|\nabla_{i_k}f(x^k)|^2.
\end{align*}
This implies a convergence rate of
\[
f(x^{k+1}) - f(x^*) \leq \left(1 - \frac{\mu_1(1-\epsilon_k)^2}{L}\right)[f(x^k) - f(x^*)].
\]
Thus, the convergence rate of the method is nearly identical to using the exact 
GS rule for small $\epsilon_k$ (and it degrades gracefully with $\epsilon_k)$. This is in 
contrast to having an error in the gradient~\citep{friedlander2011hybrid}, where 
the error $\epsilon$ must decrease to zero over time.

\subsection{Additive Errors}

In the additive error regime, the approximate GS rule chooses an $i_k$ satisfying
\[
|\nabla_{i_k} f(x^k)| \geq \norm{\nabla f(x^k)}_\infty - \epsilon_k,
\]
for some $\epsilon_k \geq 0$. In Appendix~\ref{app:additive}, we show that under this rule, we have
\begin{align*}
f(x^{k+1}) - f(x^*)  &\leq 
\left(1 - \frac{\mu_1}{L}\right)^k\left[f(x^0) - f(x^*) +A_k\right],
\end{align*}
where
\[
A_k \leq \min\left\{\sum_{i=1}^k\left(1-\frac{\mu_1}{L}\right)^{-i} \epsilon_i\frac{\sqrt{2 L_1}}{L} \sqrt{f(x^0)-f(x^*)}, \;\sum_{i = 1}^k \bigg ( 1 - \frac{\mu_1}{L} \bigg )^{-i}\left( \epsilon_i\sqrt{\frac{2}{L}}  \sqrt{f(x^0) - f(x^*)} + \frac{\epsilon_i^2}{2L} \right)\right\},
\]
where $L_1$ is the Lipschitz constant of $\nabla f$ with respect to the 1-norm. Note that $L_1$ could be substantially larger than $L$, so the second part of the maximum in $A_k$ is likely to be the smaller part unless the $\epsilon_i$ are large.
 This regime is closer to the case of having an error 
in the gradient, as to obtain convergence the $\epsilon_k$ must decrease to zero. 
This result implies that a sufficient condition for the algorithm to obtain a linear convergence 
rate is that the errors $\epsilon_k$  converge to zero at a linear rate. Further, if the 
errors satisfy $\epsilon_k = O(\rho^k)$ for some $\rho < (1-\mu_1/L)$,  then the convergence rate of the method 
is the same as if we used an exact GS rule.
On the other hand, if $\epsilon_k$ does not decrease to zero, we may end up repeatedly updating the 
same wrong coordinate and the algorithm will not converge (though we could switch to the randomized method if this is detected).

\section{Proximal-Gradient Gauss-Southwell}
\label{sec:prox}

One of the key motivations for the resurgence of interest in coordinate descent methods is 
their performance on problems of the form
\[
\min_{x\in\mathbb{R}^n} F(x) \equiv f(x) + \sum_{i=1}^n g_i(x_i),
\]
where $f$ is smooth and convex and the $g_i$ are convex, but possibly non-smooth. 
This includes problems with $\ell_1$-regularization, and optimization with lower 
and/or upper bounds on the variables. Similar to proximal-gradient methods, we can 
apply the proximal operator to the coordinate update,
$$
x^{k+1} = \prox\nolimits_{\frac{1}{L}g_{i_k}}\left[x^k - \frac{1}{L}\nabla_{i_k} f(x^k)e_{i_k}\right],
$$
where 
\[
\prox\nolimits_{\alpha g_i}[y] = \argmin{x\in\mathbb{R}^n} \frac{1}{2}\norm{x-y}^2 + \alpha g_i(x).
\]
With random coordinate selection,~\citet{richtarik2014iteration} show that this method 
has a convergence rate of
\[
\mathbb{E}[F(x^{k+1}) - F(x^*)] \leq \left(1 - \frac{\mu}{nL}\right)[F(x^k ) - F(x^*)],
\]
%where $\mu^g$ is the strong-convexity constant of the function $\sum_i g_i$. 
%We obtain~\eqref{eq:random-descent} if $\mu^g$ is zero. Although note that if $\mu^g > 0$, then 
%we can have $\mu=0$, and thus strong-convexity of $f$ is not needed. 
similar to the unconstrained/smooth case. 

There are several generalizations of the GS rule
to this scenario. Here we consider three possibilities, all of which are equivalent to the GS rule
if the $g_i$ are not present.
First, the GS-$s$ rule chooses the coordinate with the most negative directional derivative. This
strategy is popular for $\ell_1$-regularization~\citep{shevade2003simple,wu2008coordinate,li2009coordinate} and in general is given by~\citep[see][\S8.4]{bertsekas1999nonlinear}
\begin{align*}
i_k & = \argmax{i}\left\{\min_{s\in\partial g_i}|\nabla_i f(x^k) + s|\right\}.
\end{align*}
However, the length of the step ($\norm{x^{k+1}-x^k}$)
could be arbitrarily small under this choice. In contrast, the GS-$r$ rule chooses the coordinate that maximizes the length
of the step~\citep{tseng2009coordinate,dhillon2011nearest},
\begin{align*}
i_k & = \argmax{i}\left\{\left|x_i^k - \prox\nolimits_{\frac{1}{L}g_i}\left[x_i^k - \frac{1}{L}\nabla_if(x^k)\right]\right|\right\}.
\end{align*}
This rule is effective for bound-constrained problems, but it ignores the change in the non-smooth term ($g_i(x_i^{k+1})-g_i(x_k^k)$). Finally, the GS-$q$
rule maximizes progress assuming a quadratic upper bound on $f$~\citep{tseng2009coordinate},
\begin{align*}
i_k = \argmin{i}\bigg\{\min_d \big\{& f(x^k) + \nabla_if(x^k)d + \frac{L}{2}d^2 
+ g_i(x_i^k+d) - g_i(x_i^k)\big\}\bigg\}.
\end{align*}
While the least intuitive rule, the GS-$q$ rule seems to have the best theoretical properties.
Further, if we use $L_i$ in place
of $L$ in the GS-$q$ rule (which we call the GSL-$q$ strategy), then we obtain the GSL rule if the $g_i$ are not present. In contrast, using $L_i$ in place of $L$ in the GS-$r$ rule (which we call the GSL-$r$ strategy) does not
yield the GSL rule as a special case.

In Appendix~\ref{app:prox}, we show that using the GS-$q$ rule yields a convergence rate of
\begin{align*}
F(x^{k+1}) - F(x^*) &\leq  \left(1-\frac{\mu}{Ln}\right)[f(x^k)-f(x^*)],
\end{align*}
thus matching the convergence rate of randomized coordinate descent (but deterministically rather than in expectation). In contrast, in Appendix~\ref{app:prox} we also give counter-examples showing that the above rate does not hold for the GS-$s$ or the GS-$r$ rule.
Thus, any bound for the GS-$s$ or the GS-$r$ rule would be slower than the expected rate under random selection, while
the GS-$q$ rule matches this bound. It is an open problem whether the GS-$q$ rule obtains the rate $(1-\mu_1/L)$ in general, but in the next section we discuss special cases where rates depending on $\mu_1$ can be obtained.

\subsection{Rates Depending on $\mu_1$}

First, we note that if the $g_i$ are linear then the GS-$q$ rule obtains
\begin{align}
\label{eq:fastProx}
F(x^{k+1}) - F(x^*) &\leq  \left(1-\frac{\mu_1}{L}\right)[f(x^k)-f(x^*)],
\end{align}
since in this particular (smooth) case the algorithm and assumptions are identical to the setting of Section~\ref{sec:new}: the GS-q rule chooses the same coordinate to update as the GS rule applied to $F$, while $F$ has the same $L$ and $\mu_1$ as $f$ because the $g_i$ are linear.

It is possible to change the update rule in order to obtain a rate depending on $\mu_1$ for general $g$. 
In particular, after the  publication of this work~\citet{song2017accelerated} considered another generalization of the GS rule in the context of $\ell_1$-regularization. A generalized version of their update rule is given by
\begin{align*}
x^{k+1} & = x^k + d^k,\\
d^k & \in \min_{d\in\R^n}\left\{\langle \nabla f(x^k),d\rangle + \frac{L_1}{2} \norm{d}_1^2 + g(x^k + d)\right\},
\end{align*}
where $g(x) = \sum_{i=1}^n g_i(x_i)$ and where $L_1$ is the Lipschitz constant of $f$ in the $\ell_1$-norm. We call this the GS-$1$ rule, and simlar to the other GS-$*$ rules it is equivalent to the GS rule if the $g_i$ are not present. This equivalence follows from viewing the GS rule as steepest descent in the $\ell_1$-norm~\citep[][Section~9.4.2]{boyd2004convex}.~\citet[][Appendix~A.8]{nutini2018greed} shows that this rule obtains a convergence rate of
\begin{align*}
F(x^{k+1}) - F(x^*) &\leq  \left(1-\frac{\mu_1}{L_1}\right)[f(x^k)-f(x^*)].
\end{align*}
Note that $L_1 \geq L$, so this rate is slower than the other rates we have shown involving $\mu_1$. Further, unlike the other non-smooth generalizations of the GS rule, for non-linear $g$ this generalization may select more than one variable to update at each iteration. Thus it would be more appropriate to refer to this as a block coordinate descent method than a coordinate descent method. Although~\citet{song2017accelerated} give an efficient way to compute $d$ given the gradient in the case of $\ell_1$-regularization, computing $d$ for other choices of $g$ may add an additional computational cost to the method.

Finally, consider the case of $g_i$ that are piecewise-linear. Under a suitable non-degeneracy assumption, coordinate descent methods achieve a particular ``active set'' property in a finite number of iterations~\citep{wright2012,nutini2017}. Specifically, for values of $x_i^*$ that occur at non-smooth values of $g_i$, we will have $x_i^k = x_i^*$ for all sufficiently large $k$. At this point, none of the four GS-$*$ rules would select such coordinates again. Similarly, for values where $x_i^*$ occurs at smooth values of $g_i$, the iterates will eventually be confined to a region where the $g_i$ is smooth. Once this ``active set'' identification happens for piecewise-linear $g_i$, the iterates will be confined to a region where the selected $g_i$ are linear. At this point, linearity means that only one coordinate will be selected by the GS-$1$ rule and it will select the same coordinate as the GS-$q$ rule. Further, at this point the analysis of~\citet[][Appendix~A.8]{nutini2018greed} can be applied with $L$ instead $L_1$ for the GS-$q$ rule which leads to a rate of $(1-\mu_1/L)$ as in the smooth case.

\section{Experiments}

\begin{figure*}
\includegraphics[width=.49\textwidth]{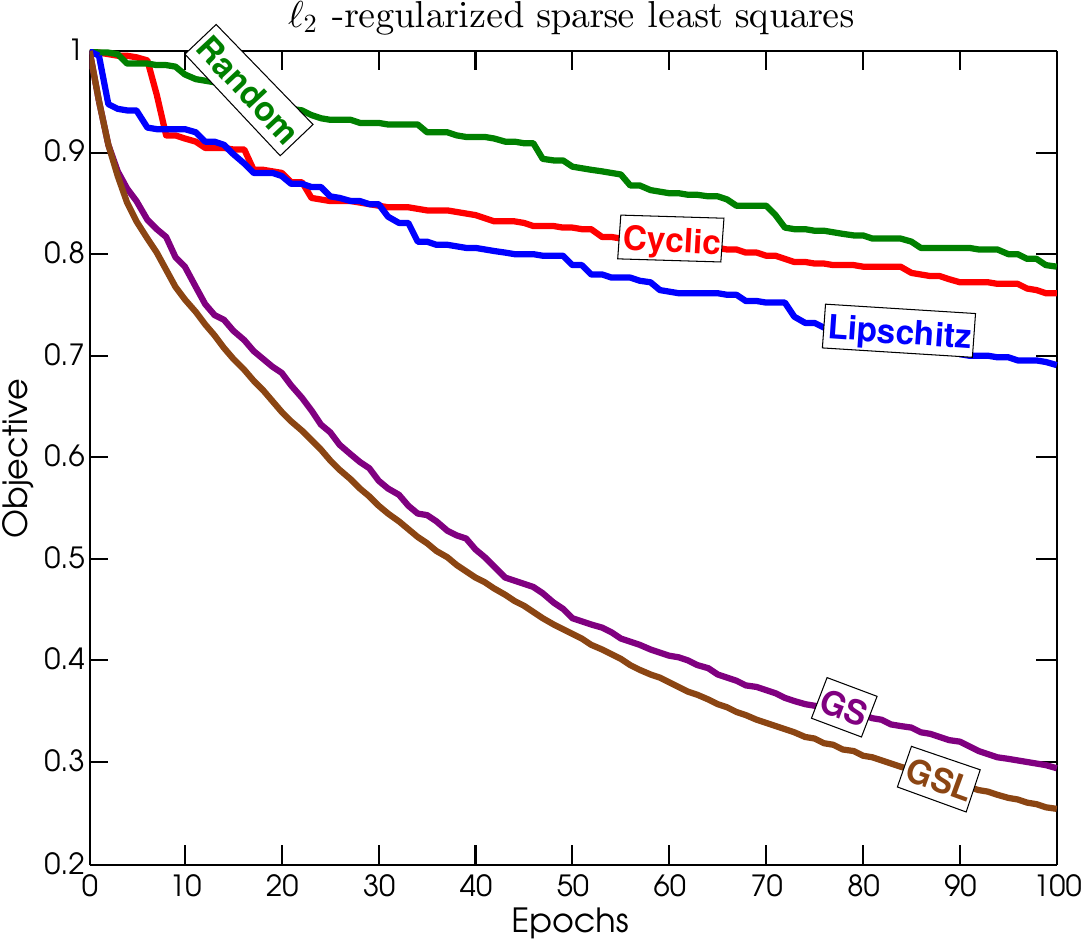}
\includegraphics[width=.49\textwidth]{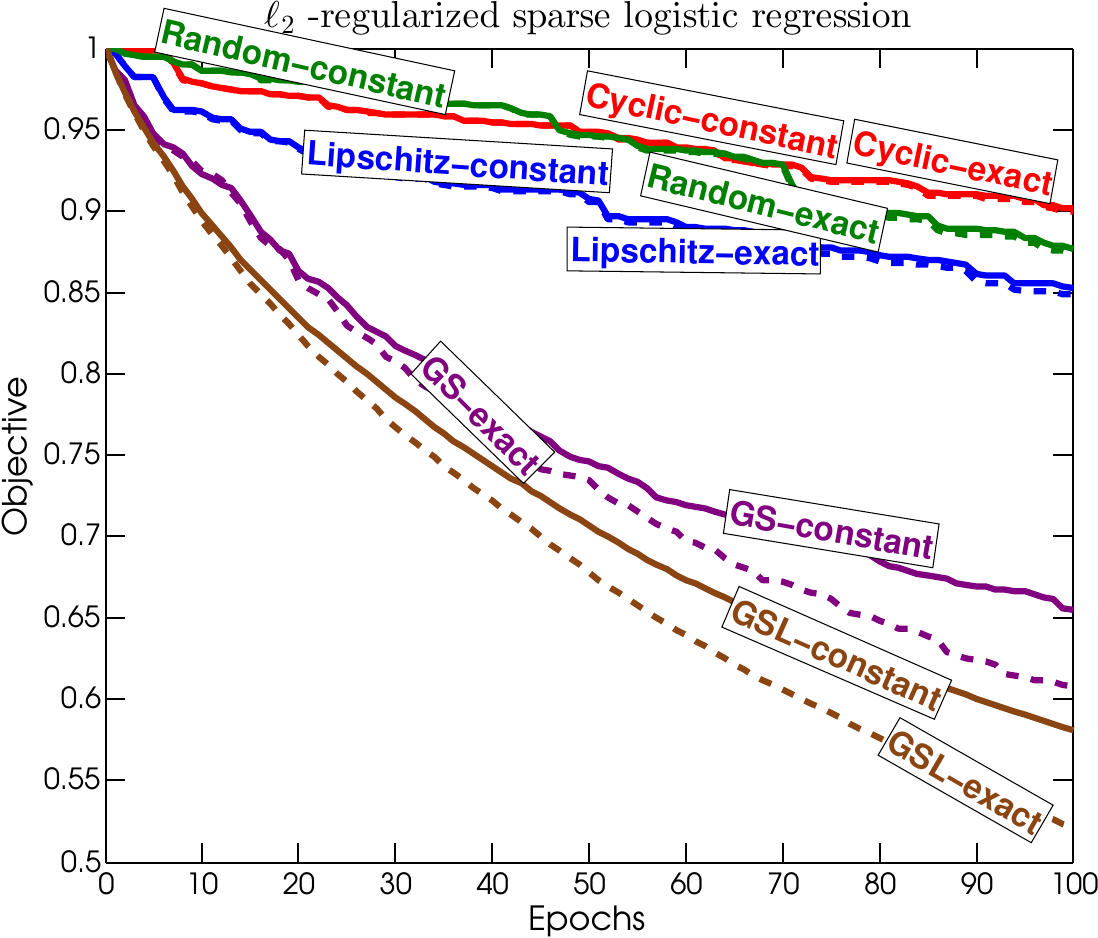}\\
\includegraphics[width=.49\textwidth]{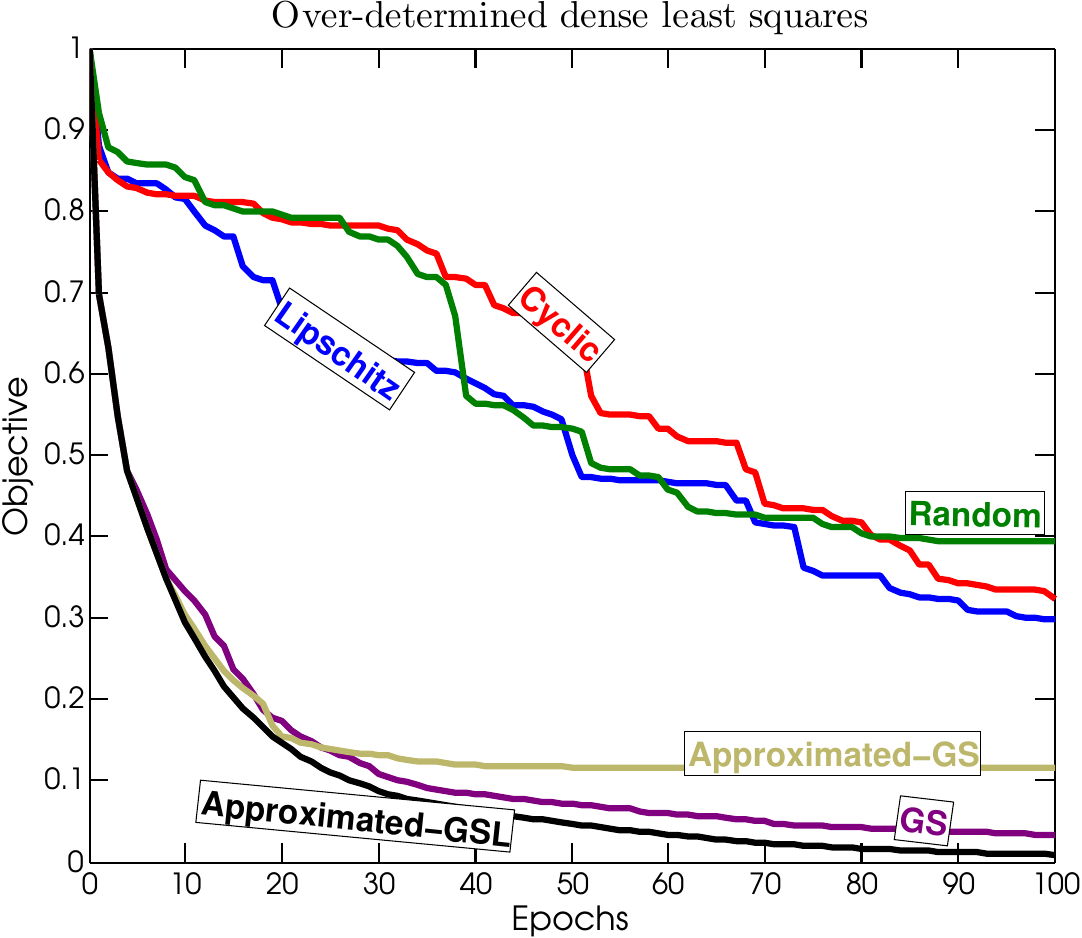}
\includegraphics[width=.49\textwidth]{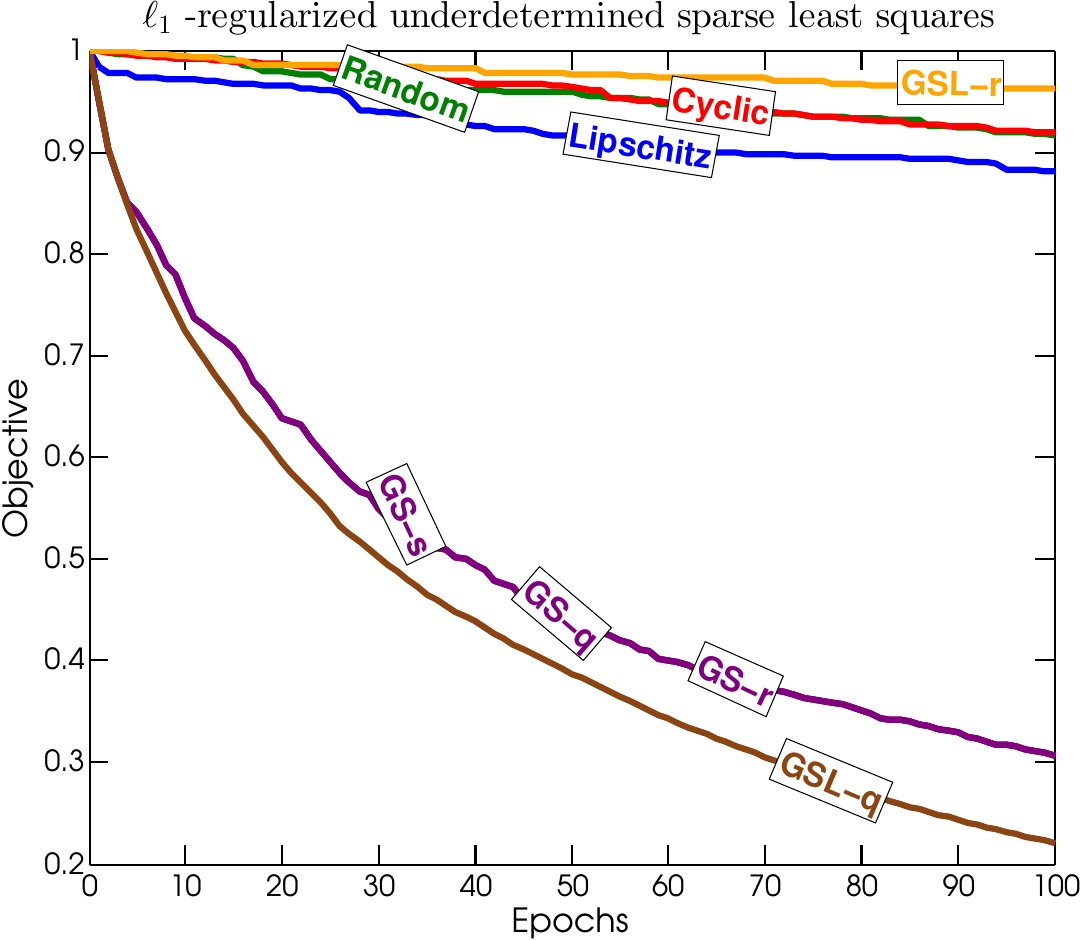}
\caption{Comparison of coordinate selection rules for 4 instances of problem $h_1$.}
\label{fig:h1}
\end{figure*}

We first compare the efficacy of different coordinate selection rules on the following simple instances of $h_1$.
\textbf{$\ell_2$-regularized sparse least squares}: Here we consider the problem
\[
\min_x \frac{1}{2m}\norm{Ax - b}^2 + \frac{\lambda}{2}\norm{x}^2,
\]
an instance of problem $h_1$. We set $A$ to be an $m$ by $n$ matrix with entries sampled from a $\mathcal{N}(0,1)$ distribution (with $m=1000$ and $n=1000$). We then added 1 to each entry (to induce a dependency between columns), multiplied each column by a sample from $\mathcal{N}(0,1)$ multiplied by ten (to induce different Lipschitz constants across the coordinates), and only kept each entry of $A$ non-zero with probability $10\log(n)/n$ (a sparsity level that allows the Gauss-Southwell rule to be applied with cost $O(\log^3(n))$. We set $\lambda=1$ and $b = Ax + e$, where the entries of $x$ and $e$ were drawn from a $\mathcal{N}(0,1)$ distribution. In this setting, we used a step-size of $1/L_i$ for each coordinate $i$, which corresponds to exact coordinate optimization.

\textbf{$\ell_2$-regularized sparse logistic regression}: Here we consider the problem
\[
\min_x \frac{1}{m}\sum_{i=1}^m\log(1+\exp(-b_ia_i^Tx)) + \frac{\lambda}{2}\norm{x}^2.
\]
We set the $a_i^T$ to be the rows of $A$ from the previous problem, and set $b = $ sign$(Ax)$, but randomly flipping each $b_i$ with probability $0.1$. In this setting, we compared using a step-size of $1/L_i$ to using exact coordinate optimization.

\textbf{Over-determined dense least squares}: Here we consider the problem
\[
\min_x \frac{1}{2m}\norm{Ax - b}^2,
\]
but, unlike the previous case, we do not set elements of $A$ to zero and we make $A$ have dimension $1000$ by $100$. Because the system is over-determined, it does not need an explicit strongly-convex regularizer to induce global strong-convexity.
In this case, the density level means that the exact GS rule is not efficient. Hence, we use a balltree structure~\citep{omohundro1989five} to implement an efficient approximate GS rule based on the connection to the NNS problem discovered by~\citet{dhillon2011nearest}. On the other hand, we can compute the exact GSL rule for this problem as a NNS problem as discussed in Section~\ref{sec:nns}.

\textbf{$\ell_1$-regularized underdetermined sparse least squares}: Here we consider the non-smooth problem
\[
\min_x \frac{1}{2m}\norm{Ax - b}^2 + \lambda\norm{x}_1.
\]
We generate $A$ as we did for the $\ell_2$-regularized sparse least squares problem, except with the dimension $1000$ by $10000$. This problem is not globally strongly-convex, but will be strongly-convex along the dimensions that are non-zero in the optimal solution.

We plot the objective function (divided by its initial value) of coordinate descent under different selection rules in Figure~\ref{fig:h1}. Even on these simple datasets, we see dramatic differences in performance between the different strategies. In particular, the GS rule outperforms random coordinate selection (as well as cyclic selection) by a substantial margin in all cases. The Lipschitz sampling strategy can narrow this gap, but it remains large (even when an approximate GS rule is used). The difference between GS and randomized selection seems to be most dramatic for the $\ell_1$-regularized problem; the GS rules tend to focus on the non-zero variables while most randomized/cyclic updates focus on the zero variables, which tend not to move away from zero.\footnote{To reduce the cost of the GS-$s$ method in this context,~\citet{shevade2003simple} consider a variant where we first compute the GS-$s$ rule for the non-zero variables and if an element is sufficiently large then they do not consider the zero variables.}
Exact coordinate optimization and using the GSL rule seem to  give modest but consistent improvements. The three non-smooth GS-$*$ rules had nearly identical performance despite their different theoretical properties. The GSL-$q$ rule gave better performance than the GS-$*$ rules, while the  the GSL-$r$ variant performed worse than even cyclic and random strategies. We found it was also possible to make the GS-$s$ rule perform poorly by perturbing the initialization away from zero. While these experiments plot the performance in terms of the number of iterations, in Appendix~\ref{app:RT} we show that the GS-$*$ rules can also be advantageous in terms of runtime.

\begin{figure}
\centering
\includegraphics[width=.5\textwidth]{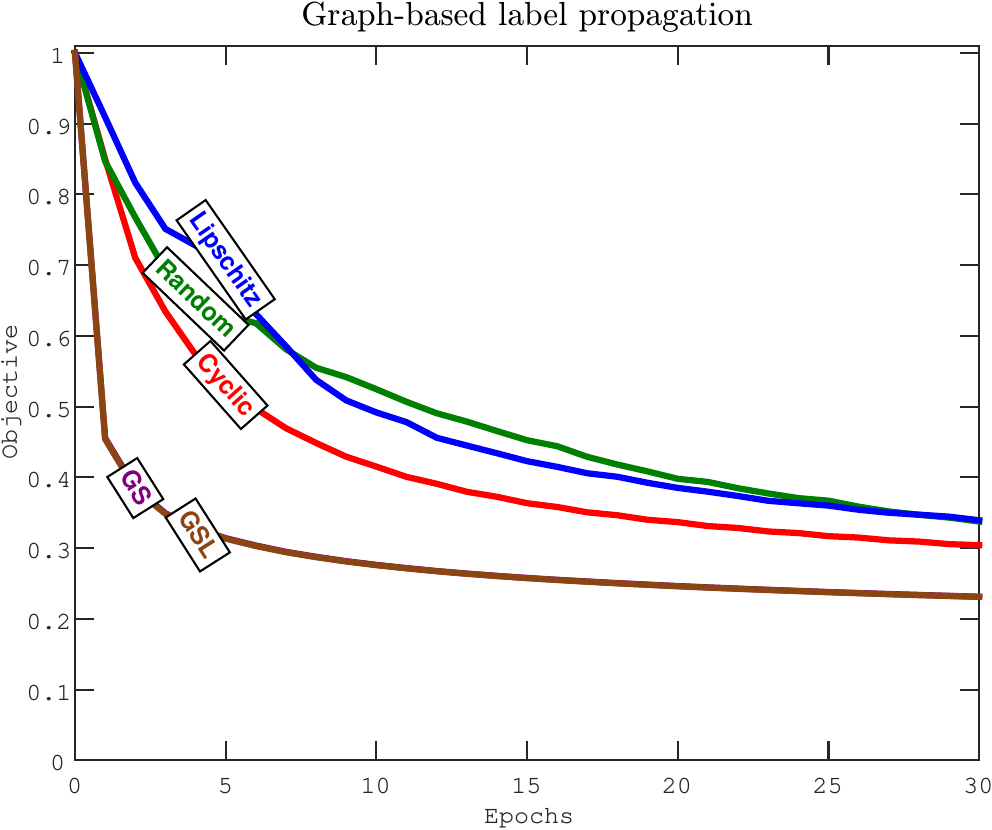}
\caption{Comparison of coordinate selection rules for graph-based semi-supervised learning.}
\label{fig:h2}
\end{figure}

We next consider an instance of problem $h_2$, performing label propagation for semi-supervised learning in the `two moons' dataset~\citep{zhou2004learning}. We generate $500$ samples from this dataset, randomly label five points in the data, and connect each node to its five nearest neighbours. This high level of sparsity is typical of graph-based methods for semi-supervised learning, and allows the exact Gauss-Southwell rule to be implemented efficiently. We use the quadratic labeling criterion of~\citet{bengio2006label}, which allows exact coordinate optimization and is normally optimized with cyclic coordinate descent. We plot the performance under different selection rules in Figure~\ref{fig:h2}. Here, we see that even cyclic coordinate descent outperforms randomized coordinate descent, but that the GS and GSL rules give even better performance. We note that the GS and GSL rules perform similarly on this problem since the Lipschitz constants do not vary much.

\section{Discussion}

It is clear that the GS rule is not practical for every problem where randomized 
methods are applicable. Nevertheless, we have shown that even approximate GS
rules can obtain better convergence rate bounds than fully-randomized methods.
 We have given a similar justification
for the use of exact coordinate optimization, and we note that our argument
could also be used to justify the use of exact coordinate optimization within
randomized coordinate descent methods (as used in our experiments). We have
also proposed the improved 
GSL rule, and considered approximate/proximal variants. We expect our analysis also applies to block updates by using mixed norms $\norm{\cdot}_{p,q}$, and could be used for accelerated/parallel methods~\citep{fercoq2013accelerated},
for primal-dual rates of dual coordinate ascent~\citep{schwartz12}, for successive projection methods~\citep{leventhal2010randomized}, for boosting algorithms~\citep{ratsch2001convergence}, and for scenarios without strong-convexity under general error bounds~\citep{luo1993error}.

\section*{Acknowledgements}

We would like to thank the anonymous referees for their useful comments that significantly improved the paper.
Julie Nutini is funded by an NSERC Canada Graduate Scholarship.

%%%%%%%%%%%%%%%%%%%%%%%%%%%%%%%%%%%%%%%%%%%%%%%%%%%

\begin{appendices}

\section{Efficient calculation of GS rules for sparse problems}
\label{app:h12}

We first give additional details on how to calculate the GS rule efficiently for sparse instances of problems $h_1$ and $h_2$. We will consider the case where each $g_i$ is smooth, but the ideas can be extended to allow a non-smooth $g_i$. Further, note that the efficient calculation does not rely on convexity, so these strategies can also be used for non-convex problems.

\subsection{Problem $h_2$}

Problem $h_2$ has the form
\[
h_2(x) := \sum_{i \in V}g_i(x_i) + \sum_{(i,j)\in E}f_{ij}(x_i,x_j),
\]
where each $g_i$ and $f_{ij}$ are differentiable and $G=\{V,E\}$ is a graph where the number of vertices $|V|$ is the same as the number of variables $n$. If all nodes in the  graph have a degree (number of neighbours) bounded above by some constant $d$, we can implement the GS rule in $O(d\log n)$ after an $O(n + |E|)$ time initialization by maintaining the following information about $x^k$:
\begin{enumerate}
\item A vector containing the values $\nabla_i g_i(x_i^k)$.
\item A matrix containing the values $\nabla_i f_{ij}(x_i^k,x_j^k)$ in the first column and $\nabla_j f_{ij}(x_i^k,x_j^k)$ in the second column.
\item The elements of the gradient vector $\nabla h_2(x^k)$ stored in a binary max heap data structure~\citep[see][Chapter~6]{cormen2001introduction}. 
\end{enumerate}
Given the heap structure, we can compute the GS rule in $O(1)$ by simply reading the index value of the root node in the max heap. The costs for initializing these structures are:
\begin{enumerate}
\item $O(n)$ to compute $g_i(x_i^0)$ for all $n$ nodes.
\item $O(|E|)$ to compute $\nabla_{ij} f_{ij}(x_i^0,x_j^0)$ for all $|E|$ edges.
\item $O(n + |E|)$ to sum the values in the above structures to compute $\nabla h(x^0)$, and $O(n)$ to construct the initial max heap.
\end{enumerate}
Thus, the one-time initialization cost is $O(n + |E|)$. The costs of updating the data structures after we update $x_{i_k}^k$ to $x_{i_k}^{k+1}$ for the selected coordinate $i_k$ are:
\begin{enumerate}
\item $O(1)$ to compute $g_{i_k}(x_{i_k}^{k+1})$.
\item $O(d)$ to compute $\nabla_{ij}f_{ij}(x_i^{k+1},x_j^{k+1})$ for $(i,j) \in E$ and $i=i_k$ or $j=i_k$ (only $d$ such values exist by assumption, and all other $\nabla_{ij}f_{ij}(x_i,x_j)$ are unchanged).
\item $O(d)$ to update up to $d$ elements of $\nabla h(x^{k+1})$ that differ from $\nabla h(x^k)$ by using differences in changed values of $g_i$ and $f_{ij}$, followed by $O(d \log n)$ to perform $d$ updates of the heap at a cost of $O(\log n)$ for each update. 
\end{enumerate}
The most expensive part of the update is modifying the heap, and thus the total cost is $O(d \log n)$.\footnote{For less-sparse problems where $n < d\log n$, using a heap is actually inefficient and we should simply store $\nabla h(x^k)$ as a vector. The initialization cost is the same, but we can then perform the GS rule in $O(n)$ by simply searching through the vector for the maximum element.}

\subsection{Problem $h_1$}

Problem $h_1$ has the form
\[
h_1(x)  := \sum_{i=1}^n g_i(x_i) + f(Ax),
\]
where $g_i$ and $f$ are differentiable, and $A$ is an $m$ by $n$ matrix where we denote column $i$ by $a_i$ and row $j$ by $a_j^T$.
Note that $f$ is a function from $\R^m$ to $\R$, and we assume $\nabla_j f$ only depends on $a_j^Tx$. While this is a strong assumption (e.g., it rules out $f$ being the product function), this class includes a variety of notable problems like the least squares and logistic regression models from our experiments. If $A$ has $z$ non-zero elements, with a maximum of $c$ non-zero elements in each column and $r$ non-zero elements in each row, then with a pre-processing cost of $O(z)$ we can implement the GS rule in this setting in $O(cr\log n)$ by maintaining the following information about $x^k$:
\begin{enumerate}
\item A vector containing the values $\nabla_i g_i(x_i^k)$.
\item A vector containing the product $Ax^k$.
\item A vector containing the values $\nabla f(Ax^k)$.
\item A vector containing the product $A^T\nabla f(Ax^k)$.
\item The elements of the gradient vector $\nabla h_1(x^k)$ stored in a binary max heap data structure.
\end{enumerate}
The heap structure again allows us to compute the GS rule in $O(1)$, and the costs of initializing these structures are:
\begin{enumerate}
\item $O(n)$ to compute $g_i(x_i^0)$ for all $n$ variables.
\item $O(z)$ to compute the product $Ax^0$.
\item $O(m)$ to compute $\nabla f(Ax^0)$ (using that $\nabla_j f$ only depends on $a_j^Tx^0$).
\item $O(z)$ to compute $A^T\nabla f(Ax^0)$.
\item $O(n)$ to add the $\nabla_i g_i(x_i^0)$ to the above product to obtain $\nabla h_1(x^0)$ and construct the initial max heap.
\end{enumerate}
As it is reasonable to assume that $z \geq m$ and $z \geq n$ (e.g., we have at least one non-zero in each row and column), the cost of the initialization is thus $O(z)$. The costs of updating the data structures after we update $x_{i_k}^k$ to $x_{i_k}^{k+1}$ for the selected coordinate $i_k$ are:
\begin{enumerate}
\item $O(1)$ to compute $g_{i_k}(x_{i_k}^{k+1})$.
\item $O(c)$ to update the product using $Ax^{k+1} = Ax^k + (x_{i_k}^{k+1}-x_{i_k}^k)a_i$, since $a_i$ has at most $c$ non-zero values.
\item $O(c)$ to update up to $c$ elements of $\nabla f(Ax^{k+1})$ that have changed (again using that $\nabla_j f$ only depends on $a_j^Tx^{k+1}$). 
\item $O(cr)$ to perform up to $c$ updates of the form $A^T\nabla f(Ax^{k+1}) = A^T\nabla f(Ax^k) + (\nabla_j f(Ax^{k+1})- \nabla_j f(Ax^k))(a_i)^T$, where each update costs $O(r)$ since each $a_i$ has at most $r$ non-zero values.
\item $O(cr\log n)$ to update the gradients in the heap.
\end{enumerate}
The most expensive part is again the heap update, and thus the total cost is $O(cr\log n)$.
%\footnote{If $n < cr\log n$, we do not need a heap as we can simply store the gradient and search through it to reduce the cost to $O(n + cr)$.}

\section{Relationship between $\mu_1$ and $\mu$}
\label{app:mu12}

We can establish the relationship between $\mu$ and $\mu_1$ by using the known relationship between the $2$-norm and the $1$-norm,
\[
\norm{x}_1 \geq \norm{x} \geq \frac{1}{\sqrt{n}}\norm{x}_1.
\]
In particular, if we assume that $f$ is $\mu$-strongly convex in the $2$-norm, then for all $x$ and $y$ we have
\begin{align*}
f(y) & \geq f(x) + \langle \nabla f(x),y-x\rangle + \frac{\mu}{2}\norm{y-x}^2\\
& \geq f(x) + \langle \nabla f(x),y-x\rangle + \frac{\mu}{2n}\norm{y-x}_1^2,
\end{align*}
implying that $f$ is at least $\frac{\mu}{n}$-strongly convex in the $1$-norm. Similarly, if we assume that a given $f$ is $\mu_1$-strongly convex in the $1$-norm then for all $x$ and $y$ we have
\begin{align*}
f(y) & \geq f(x) + \langle \nabla f(x),y-x\rangle + \frac{\mu_1}{2}\norm{y-x}_1^2\\
& \geq f(x) + \langle \nabla f(x),y-x\rangle + \frac{\mu_1}{2}\norm{y-x}^2,
\end{align*}
implying that $f$ is at least $\mu_1$-strongly convex in the $2$-norm. Summarizing these two relationships, we have
\[
\frac{\mu}{n} \leq \mu_1 \leq \mu.
\]

\section{Analysis for separable quadratic case}
\label{app:mu1}

We first establish an equivalent definition of strong-convexity in the $1$-norm, along the lines of~\citet[Theorem 2.1.9]{Nes04b}. Subsequently, we use this equivalent definition to derive $\mu_1$ for a separable quadratic function.

\subsection{Equivalent definition of strong-convexity}

Assume that $f$ is $\mu_1$-strongly convex in the $1$-norm, so that for any $x,y \in \R^n$  we have
\[
	f(y) \ge f(x) + \langle \nabla f(x), y - x \rangle + \frac{\mu_1}{2} \| y - x \|^2_1.
\]
Reversing $x$ and $y$ in the above gives
\[
	f(x) \ge f(y) + \langle \nabla f(y), x - y \rangle + \frac{\mu_1}{2} \| x - y \|^2_1,
\]
and adding these two together yields
\begin{equation}
\label{eq:sc2}
	\langle \nabla f(y) - \nabla f(x), y - x \rangle \ge \mu_1 \| y - x \|^2_1.
\end{equation}
Conversely, assume  that for all $x$ and $y$ we have
\[
	\langle \nabla f(y) - \nabla f(x), y - x \rangle \ge \mu_1 \| y - x \|^2_1,
\]
and consider the function $g(\tau) = f(x+\tau(y - x))$ for $\tau \in \R$. Then
\begin{align*}
	f(y) - f(x) - \langle \nabla f(x), y - x \rangle
	&= g(1) - g(0) - \langle \nabla f(x), y - x \rangle\\
	&= \int_0^1 \frac{dg}{d\tau}(\tau) - \langle \nabla f(x), y - x \rangle ~d\tau \\
	&= \int_0^1 \langle \nabla f(x + \tau(y - x)), y - x \rangle - \langle \nabla f(x), y - x \rangle ~d\tau \\
	&= \int_0^1 \langle \nabla f(x + \tau(y - x)) - \nabla f(x), y - x \rangle ~d\tau \\
	&\ge \int_0^1 \frac{\mu_1}{\tau} \| \tau(y - x) \|^2_1 ~d\tau \\
	&= \int_0^1 \mu_1 \tau \| y - x \|^2_1 ~d\tau \\
	&= \frac{\mu_1}{2} \tau^2 \| y - x \|^2_1 \bigg |^1_0 \\
	&= \frac{\mu_1}{2} \| y - x \|^2_1.
\end{align*}
Thus, $\mu_1$-strong convexity in the $1$-norm is equivalent to having
\begin{equation}\label{first}
	\langle \nabla f(y) - \nabla f(x), y - x \rangle \ge \mu_1 \| y - x \|^2_1 \quad \forall ~x,y.
\end{equation}

\subsection{Strong-convexity constant $\mu_1$ for separable quadratic functions}

%Let $\mu_1 > 0$ be the strong convexity constant of $f$ with respect to the $1$-norm $\| \cdot \|_1$, 
%i.e., for any $x, y \in \R^n$
%$$
%	f(y) \ge f(x) + \langle \nabla f(x), y - x \rangle + \frac{\mu_1}{2} \| y - x \|^2_1.
%$$

%\begin{proof}

Consider a strongly convex quadratic function $f$ with a diagonal Hessian 
$H = \nabla^2 f(x) = \diag{\lambda_1, \dots, \lambda_n}$, where $\lambda_i > 0$ for all $i = 1, \dots, n$. We show that in this case
$$
	\mu_1 = \left(\sum_{i=1}^n\frac{1}{\lambda_i}\right)^{-1}.
$$ 
From the previous section, $\mu_1$ is the minimum value such that~\eqref{first} holds,
\[
\mu_1 = \inf_{x\neq y} \frac{\langle \nabla f(y) - \nabla f(x), y - x \rangle}{\norm{y-x}_1^2}.
\] 
Using $\nabla f(x) = Hx + b$ for some $b$ and letting $z = y-x$, we get
\begin{align*}
\mu_1 & = \inf_{x \neq y} \frac{\langle (Hy - b) - (Hx- b),y-x\rangle}{\norm{y-x}_1^2} \\
& = \inf_{x\neq y} \frac{\langle H(y-x),y-x\rangle}{\norm{y-x}_1^2}\\
& = \inf_{z\neq 0} \frac{z^THz}{\norm{z}_1^2}\\
	&= \min_{\|z\|_1 = 1} z^T H z \\
	&= \min_{e^Tz = 1} \sum_{i=1}^n \lambda_i z_i^2, 
\end{align*}
where the last two lines use that the objective is invariant to scaling of $z$ and to the sign of $z$ (respectively), and where $e$ is a vector containing a one in every position. This is an equality-constrained strictly-convex quadratic program, so its solution is given as a stationary point $(z^*,\eta^*)$ of the Lagrangian,\[
	\Lambda(z,\eta) = \sum_{i = 1}^n \lambda_i z_i^2 + \eta (1- e^Tz).
\]
Differentiating with respect to each $z_i$ for $i = 1, \dots, n$ and equating to zero, we have for all $i$ that $2 \lambda_i z_i^* - \eta^*= 0$, or
\begin{equation}\label{eta}
z_i^* = \frac{\eta^*}{2\lambda_i}.
\end{equation}
Differentiating the Lagrangian with respect to $\eta$ and equating to zero we obtain $1- e^Tz^* = 0$, or equivalently
\[
	1 = e^T z^* = \frac{\eta^*}{2}\sum_j \frac{1}{\lambda_j},
\]
which yields 
\[
	\eta^* = 2 \left(\sum_j \frac{1}{\lambda_j}\right)^{-1}.
\]
Combining this result for $\eta^*$ with equation \eqref{eta}, we have
\[
	z_i^* = \frac{1}{\lambda_i} \left(\sum_j \frac{1}{\lambda_j}\right)^{-1}.
\]
This gives the minimizer, so we evaluate the objective at this point to obtain $\mu_1$,
\begin{align*}
	\mu_1 
	& = \sum_{i=1}^n \lambda_i (z_i^*)^2\\
	&= \sum_{i=1}^n \lambda_i \left( \frac{1}{\lambda_i} \left(\sum_{j=1}^n \frac{1}{\lambda_j}\right)^{-1}\right)^2\\
	&= \sum_{i=1}^n \frac{1}{\lambda_i}\left(\sum_{j=1}^n \frac{1}{\lambda_j}\right)^{-2}\\
	&= \left(\sum_{j=1}^n \frac{1}{\lambda_j}\right)^{-2}\left(\sum_{i=1}^n \frac{1}{\lambda_i}\right)\\
	& = \left(\sum_{j=1}^n \frac{1}{\lambda_j}\right)^{-1}.
\end{align*}
%which is equivalent to solving
%\begin{equation*}
%\begin{split}
%	\min & \quad y \\
%	\text{s.t.}   & \quad - y\e \le Hz \le y\e, \\
%			& \quad \| z \|_1 = 1,
%\end{split}
%\end{equation*}
%where $y \in \R$ and $\e = [1, \dots, 1]^T \in \R^n$. From the first inequality, we get
%$$
%	- yH^{-1}\e \le z \le yH^{-1}\e \quad \Rightarrow \quad | z_i | \le y\frac{1}{\lambda_i}, \quad \forall i.
%$$
%Given that we are minimizing $y$, and combining this with the equality constraint on $z$, we have
%\begin{equation}\label{ybound:1}
%	\|z \|_1 
%	= y \bigg [ \frac{1}{\lambda_1} + \dots + \frac{1}{\lambda_n} \bigg] 
%	= 1 \quad \Rightarrow \quad y = \frac{1}{\frac{1}{\lambda_1} + \dots + \frac{1}{\lambda_n}}.
%\end{equation}
%Thus,
%$$
%	\mu_1 
%	= \frac{1}{\frac{1}{\lambda_1} + \dots + \frac{1}{\lambda_n}} 
%	= \frac{ \prod_{i = 1}^n \lambda_i}{ \sum_{j = 1}^n \prod_{i \not = j} \lambda_i}.
%$$
%\end{proof}

%%%%%%%%%%%%%%%%%%%%%%%%%%%%%%%%%%%%%%%%%%%%%%%%%%%

%\subsection*{Behaviour of Algorithms for Separable Quadratic}

%Consider the case where $\lambda_1 = \beta$ and $\lambda_2 = \lambda_3 = \cdots \lambda_n = \alpha$. In this case, whenever we select to update coordinate $i$ for $i > 0$, using a step-size of $1/L= 1/\alpha$ the algorithm moves variables $x_i$ to its optimal value. In contrast, when we select to update coordinate $1$, the step-size of $1/\alpha$ is too small.

%\subsection*{Example showing advantage of Gauss-Southwell}

\section{Gauss-Southwell with exact optimization}
\label{app:exact}

We can obtain a faster convergence for GS using exact coordinate optimization for sparse variants of problems $h_1$ and $h_2$, by observing that the convergence rate can be expressed in terms of the sequence of $(1-\mu_1/L_{i_k})$ values,
\[
f(x^k) - f(x^*) \leq \left[\prod_{j=1}^k \left(1-\frac{\mu_1}{L_{i_j}}\right)\right][f(x^0) - f(x^*)].
\]
The worst case occurs when the product of the $(1-\mu_1/L_{i_k})$ values is as large as possible. However, using exact coordinate optimization guarantees that, after we have updated coordinate $i$, the GS rule will never select it again until one of its neighbours has been selected. Thus, we can obtain 
a tighter bound on the worst-case convergence rate using GS with exact 
coordinate optimization on iteration $k$, by solving the following 
combinatorial optimization problem defined on a weighted graph:
\begin{problem}
We are given a graph $G=(V,E)$ with $n$ nodes, a number $M_i$ associated with 
each node $i$, and an iteration number $k$. Choose a sequence $\{i_t\}_{t=1}^k$ 
that maximizes the sum of the $M_{i_t}$, subject to the following constraint: after each time node 
$i$ has been chosen, it cannot be chosen again until after a neighbour of node $i$ 
has been chosen.
\end{problem}
\noindent
We can use the $M_i$ chosen by this problem to obtain an upper-bound on the sequence of $\log(1-\mu_1/L_i)$ values, and if the largest $M_i$ values are not close to each 
other in the graph, then this rate can be much faster than the rate obtained by 
alternating between the largest $M_i$ values. In the particular case of chain-structured graphs, a worst-case sequence can be constructed that spends all but $O(n)$ iterations in one of two solution modes: (i) alternate between two nodes $i$ and $j$ that are
connected by an edge with the highest value of $\frac{M_i+M_j}{2}$, or (ii) alternate between 
three nodes $\{i,j,k\}$ with the highest value of $\frac{M_i+M_j+M_k}{3}$, where there is an edge from 
$i$ to $j$ and from $j$ to $k$, but not from $i$ to $k$. To show that these are the 
two solution modes, observe that the solution must eventually cycle because there are a finite number of nodes. If you have more 
than three nodes in the cycle, then you can always remove one node from the cycle 
to obtain a better average weight for the cycle without violating the constraint. 
%However, the solution cannot be a 
%3-clique, because otherwise you could remove the lowest-weight node in the clique 
%and obtain a higher-average-weight cycle without violating the constraint. 
We will fall into mode (i) if the average of 
$M_i$ and $M_j$ in this mode is larger than the average of $M_i$, $M_j$ and $M_k$ 
in the second mode. We can construct a solution to this problem that consists of a 
`burn-in' period, where we choose the largest $M_i$, followed by repeatedly going 
through the better of the two solution modes up until the final three steps, where a 
`burn-out' phase arranges to finish with several large $M_i$. By setting $M_i = \log(1-\mu_1/L_i)$, this leads to a convergence
rate of the form
\[
 f(x^k) - f(x^*) \leq O\left(\max\{\rho_2^G,\rho_3^G\}^k\right)[f(x^0) - f(x^*)],
\]
where $\rho_2^G$ is the maximizer of $\sqrt{(1-\mu_1/L_i)(1-\mu_1/L_j)}$ among all consecutive nodes $i$ and $j$ in the chain, and $\rho_3^G$ is the maximizer of $\sqrt[3]{(1-\mu_1/L_i)(1-\mu_1/L_j)(1-\mu_1/L_k)}$ among consecutive nodes $i$, $j$, and $k$.
%where $L_2^G$ is the maximum average Lipschitz constant among neighbours in the graph, and $L_3^G$
%is the largest average Lipschitz constant among chains of three nodes in the graph.
 The $O()$
notation gives the constant due to choosing higher $(1-\mu_1/L_i)$ values during the burn-in and burn-out periods. %, and the fact that we may end the solution mode with an $(1-\mu/L_i)$ value that is larger than $\rho_2^G$ and $\rho_3^G$.
The implication 
of this result is that, if the large $L_i$ values are more than two edges away from each other in the 
graph, then the convergence rate can be much faster.

%%%%%%%%%%%%%%%%%%%%%%%%%%%%%%%%%%%%%%%%%%%%%%%%%%%

\section{Gauss-Southwell-Lipschitz rule: convergence rate}
\label{app:GSL}

The coordinate-descent method with a constant step-size of $L_{i_k}$ uses the iteration
\[
x^{k+1} = x^k - \frac{1}{L_{i_k}}\nabla_{i_k} f(x^k)e_{i_k}.
\]
Because $f$ is coordinate-wise $L_{i_k}$-Lipschitz continuous, we obtain the following bound on the progress 
made by each iteration:
\begin{equation} \label{eq:descent_L}
\begin{aligned}
f(x^{k+1}) & \leq f(x^k) + \nabla_{i_k}f(x^k)(x^{k+1}-x^k)_{i_k} + \frac{L_{i_k}}{2}(x^{k+1} - x^k)_{i_k}^2\\
& = f(x^k) -\frac{1}{L_{i_k}} (\nabla_{i_k}f(x^k))^2 + \frac{L_{i_k}}{2}\left[\frac{1}{L_{i_k}}\nabla_{i_k}f(x^k)\right]^2 \\
& = f(x^k) - \frac{1}{2L_{i_k}}[\nabla_{i_k}f(x^k)]^2\\
& = f(x^k) - \frac{1}{2} \bigg [\frac{\nabla_{i_k}f(x^k)}{\sqrt{L_{i_k}}} \bigg]^2.
\end{aligned}
\end{equation}
By choosing the coordinate to update according to the Gauss-Southwell-Lipchitz (GSL) rule, 
\[
i_k = \argmax{i} \frac{|\nabla_i f(x^k)|}{\sqrt{L_i}},
\]
we obtain the tightest possible bound on~\eqref{eq:descent_L}. We define the following norm,
\begin{equation}
\label{eq:Lnorm}
\norm{x}_L = \sum_{i = 1}^n \sqrt{L_i}|x_i|, 
\end{equation}
which has a dual norm of
\[
\norm{x}_L^* =\max_{i} \frac{1}{\sqrt{L_i}}|x_i|.
\]
Under this notation, and using the GSL rule,~\eqref{eq:descent_L} becomes
\[
f(x^{k+1}) \leq f(x^k) - \frac{1}{2} \big (\norm{\nabla f(x^k)}^*_L \big )^2,
\]
Measuring strong-convexity in the norm $\norm{\cdot}_L$ we get
\[
f(y) \geq f(x) + \langle \nabla f(x),y-x\rangle + \frac{\mu_L}{2}\norm{y-x}_L^2.
\]
 Minimizing both sides with respect to $y$ we get
\begin{align*}
f(x^*) & \geq f(x) - \sup_y\{\langle -\nabla f(x),y-x\rangle -  \frac{\mu_L}{2}\norm{y-x}_L^2\}\\
& = f(x) - \left(\frac{\mu_L}{2}\norm{\cdot}_L^2\right)^*(-\nabla f(x)) \\
& = f(x) - \frac{1}{2\mu_L}\big (\norm{\nabla f(x)}^*_L \big)^2.
\end{align*}
Putting these together yields
\begin{equation} \label{eq:gs-new_L}
f(x^{k+1}) - f(x^*) \leq (1 - \mu_L)[f(x^k) - f(x^*)].
\end{equation}
%If all $L$ are equal, then $\mu_L = \mu_1/L_i$, so the rate is the same as before. Otherwise 
%this rate should be faster.

%%%%%%%%%%%%%%%%%%%%%%%%%%%%%%%%%%%%%%%%%%%%%%%%%%%

\section{Comparing $\mu_L$ to $\mu_1$ and $\mu$}
\label{app:muL}

By the logic Appendix \ref{app:mu12}, to establish a relationship between different strong-convexity constants under different norms, it is sufficient to establish the relationships between the squared norms. In this section, we use this to establish the relationship between $\mu_L$ defined in~\eqref{eq:Lnorm} and both $\mu_1$ and $\mu$.

\subsection{Relationship between $\mu_L$ and $\mu_1$}

We have
\[
c\norm{x}_1 - \norm{x}_L = c\sum_i |x_i| - \sum_i \sqrt{L_i}|x_i| = \sum_i (c - \sqrt{L_i})|x_i|,
\]
Assuming $c \geq \sqrt{L}$, where $L = \max_i \{ L_i \}$, the expression is non-negative and we get
$$
	\norm{x}_L \leq \sqrt{L}\norm{x}_1.
	$$
By using
\[
c\norm{x}_L - \norm{x}_1 = \sum_i (c\sqrt{L_i} - 1)|x_i|,
\]
and assuming $\displaystyle c \geq \frac{1}{\sqrt{L_{min}}}$, where $L_{min} = \min_i \{ L_i \}$, this expression is nonnegative and we get
$$
	\norm{x}_1 \leq \frac{1}{\sqrt{L_{min}}}\norm{x}_L.
	$$
	The relationship between $\mu_L$ and $\mu_1$ is based on the squared norm, so in summary we have
\[
\frac{\mu_1}{L} \leq \mu_L \leq \frac{\mu_1}{L_{min}}.
\]

%%%%%%%%%
\subsection{Relationship between $\mu_L$ and $\mu$}

Let $\vec{L}$ denote a vector with elements $\sqrt{L_i}$, and we note that 
$$
	\norm{\vec{L}} = \bigg (\sum_i (\sqrt{L_i})^2 \bigg)^{1/2} = \bigg ( \sum_i L_i \bigg)^{1/2} = \sqrt{n\bar{L}}, \quad  \text{where } \bar{L} = \frac{1}{n} \sum_i L_i.
$$
\\
Using this, we have
\[
\norm{x}_L = x^T(\sign(x)\circ\vec{L}) \leq \norm{x}\norm{\sign(x)\circ \vec{L}} = \sqrt{n\bar{L}}\norm{x}.
\]
This implies that
\[
\frac{\mu}{n\bar{L}} \leq \mu_L.
\]
Note that we can also show that $\mu_L \leq \frac{\mu}{L_{min}}$, but this is less tight than the upper bound from the previous section because  $\mu_1 \leq \mu$.

%%%%%%%%%%%%%%%%%%%%%%%%%%%%%%%%%%%%%%%%%%%%%%%%%%%

\section{Approximate Gauss-Southwell with additive error}
\label{app:additive}

%Put Julie's derivation of rate under additive error here. Use something like that $\norm{\nabla f(x^k) - \nabla f(x^*)} \leq L\norm{x^k - x^*} \leq L\norm{x^0 - x^*}$ to get rid of the gradient term?
In the additive error regime, the approximate Gauss-Southwell rule chooses an $i_k$ satisfying
$$
	| \nabla_{i_k} f(x^k) | \ge \| \nabla f(x^k) \|_\infty - \epsilon_k, \quad \text{where } \epsilon_k \geq 0~~ \forall k,
$$
and we note that we can assume $\epsilon_k \leq \norm{\nabla f(x^k)}_\infty$ without loss of generality because we must always choose an $i$ with $|\nabla_{i_k}f(x^k)| \geq 0$. 
Applying this to our bound on the iteration progress, we get
\begin{equation}
\label{eq:additive}
\begin{aligned}
	f(x^{k+1}) & \le f(x^k) - \frac{1}{2L} \bigg [ \nabla_{i_k} f(x^k ) \bigg ]^2 \\
			& \le f(x^k) - \frac{1}{2L} \big (\| \nabla f(x^k) \|_\infty - \epsilon_k \big )^2 \\
			& = f(x^k) - \frac{1}{2L} \big (\| \nabla f(x^k) \|_\infty^2 - 2 \epsilon_k \| \nabla f(x^k) \|_\infty + \epsilon_{k}^2 \big ) \\
			& = f(x^k) - \frac{1}{2L} \| \nabla f(x^k) \|_\infty^2 + \frac{\epsilon_{k}}{L}  \| \nabla f(x^k) \|_\infty - \frac{\epsilon_{k}^2}{2L} \\
			%& \le f(x^k) - \frac{1}{2L} \| \nabla f(x^k) \|_\infty^2 + \frac{\epsilon_{k}}{L}  \| \nabla f(x^k) \|_\infty. \\
%			& = f(x^k) - \frac{1}{2L} \| \nabla f(x^k) \|_\infty^2 + \frac{\epsilon_{i_k}}{2L}  \bigg ( 2 \| \nabla f(x^k) \|_\infty - \epsilon_{i_k} \bigg ) . \\
%			& \le f(x^k) - \frac{1}{2L} \| \nabla f(x^k) \|_\infty^2 + \frac{\epsilon_{i_k}}{L} \| \nabla f(x^k) \|_\infty.  \\
\end{aligned}
\end{equation}
We first give a result that assumes $f$ is $L_1$-Lipschitz continuous in the $1$-norm. This implies an inequality that we prove next, followed by a convergence rate that depends on $L_1$. However, note that $L \leq L_1 \leq Ln$, so this potentially introduces a dependency on $n$. We subsequently give a slightly less concise result that has a worse dependency on $\epsilon$ but does not rely on $L_1$.

\subsection{Gradient bound in terms of $L_1$}

We say that $\nabla f$ is $L_1$-Lipschitz continuous in the $1$-norm if we have for all $x$ and $y$ that
\[
\norm{\nabla f(x) - \nabla f(y)}_\infty \leq L_1\norm{x-y}_1.
\]
Similar to~\citet[][Theorem~2.1.5]{Nes04b}, we now show that this implies
\begin{equation}
\label{eq:L1grad}
f(y) \geq f(x) + \langle \nabla f(x),y-x\rangle + \frac{1}{2L_1}\norm{\nabla f(y) - \nabla f(x)}_\infty^2,
\end{equation}
and subsequently that
\begin{equation}
\label{eq:Linf-bound}
\norm{\nabla f(x^k)}_\infty = \norm{\nabla f(x^k) - \nabla f(x^*)}_\infty \leq \sqrt{2L_1(f(x^k)-f(x^*))} \leq \sqrt{2L_1(f(x^0)-f(x^*))},
\end{equation}
where we have used that $f(x^k) \leq f(x^{k-1})$ for all $k$ and any choice of $i_{k-1}$ (this follows from the basic bound on the progress of coordinate descent methods).

We first show that $\nabla f$ being $L_1$-Lipschitz continuous in the $1$-norm implies that
\[
f(y) \leq f(x) + \langle \nabla f(x),y-x\rangle + \frac{L_1}{2}\norm{y-x}_1^2,
\]
for all $x$ and $y$.
Consider the function $g(\tau) = f(x+\tau(y - x))$ with $\tau \in \R$. Then
\begin{align*}
	f(y) - f(x) - \langle \nabla f(x), y - x \rangle
	&= g(1) - g(0) - \langle \nabla f(x), y - x \rangle\\
	&= \int_0^1 \frac{dg}{d\tau}(\tau) - \langle \nabla f(x), y - x \rangle ~d\tau \\
	&= \int_0^1 \langle \nabla f(x + \tau(y - x)), y - x \rangle - \langle \nabla f(x), y - x \rangle ~d\tau \\
	&= \int_0^1 \langle \nabla f(x + \tau(y - x)) - \nabla f(x), y - x \rangle ~d\tau \\
	&\le \int_0^1 \| \nabla f(x + \tau(y - x)) - \nabla f(x) \|_\infty \|y - x \|_1 ~d\tau \\
	&\le \int_0^1 L_1 \tau \| y - x \|_1^2 ~d\tau \\
	&= \frac{L_1}{2} \tau^2 \| y - x \|^2_1 \bigg |^1_0 \\
	&= \frac{L_1}{2} \| y - x \|^2_1.
\end{align*}
To subsequently show~\eqref{eq:L1grad}, fix $x \in \R^n$ and consider the function 
\[
	\phi(y) = f(y) - \langle \nabla f(x),y \rangle,
\]
which is convex on $\R^n$ and also has an $L_1$-Lipschitz continuous gradient in the $1$-norm, as
\begin{align*}
	\| \phi'(y) - \phi'(x) \|_\infty 
	&= \| (\nabla f(y) - \nabla f(x)) - (\nabla f(x) - \nabla f(x))  \|_\infty \\
	&= \| \nabla f(y) - \nabla f(x) \|_\infty \\
	&\le L_1 \| y - x \|_1.
\end{align*}
As the minimizer of $\phi$ is $x$ (i.e., $\phi'(x) = 0$), for any $y \in \R^n$ we have
\begin{align*}
	\phi(x)  = \min_v \phi(v) 
	&\le \min_v \phi(y) + \langle \phi'(y), v - y \rangle + \frac{L_1}{2} \| v - y \|^2_1 \\
	&=  \phi(y) - \sup_v \langle -\phi'(y), v - y \rangle - \frac{L_1}{2} \| v - y \|^2_1 \\
	&=  \phi(y) - \frac{1}{2L_1} \| \phi'(y) \|^2_\infty.
\end{align*}
Substituting in the definition of $\phi$, we have
\begin{align*}
	f(x) - \langle \nabla f(x),x \rangle &\le f(y) - \langle \nabla f(x),y \rangle - \frac{1}{2L_1} \| \nabla f(y) - \nabla f(x) \|^2_\infty \\
	\iff \hspace{2.5cm} f(x) &\le f(y) + \langle \nabla f(x),x - y \rangle - \frac{1}{2L_1} \| \nabla f(y) - \nabla f(x) \|^2_\infty \\
	\iff \hspace{2.6cm}\! f(y) &\ge f(x) + \langle \nabla f(x), y - x \rangle + \frac{1}{2L_1} \| \nabla f(y) - \nabla f(x) \|^2_\infty.
\end{align*}

\subsection{Additive error bound in terms of $L_1$}

Using~\eqref{eq:Linf-bound} in~\eqref{eq:additive} and noting that $\epsilon_k \geq 0$, we obtain
\begin{align*}
f(x^{k+1}) & \leq f(x^k) - \frac{1}{2L} \| \nabla f(x^k) \|_\infty^2 + \frac{\epsilon_{k}}{L}  \| \nabla f(x^k) \|_\infty - \frac{\epsilon_{k}^2}{2L} \\
& \leq f(x^k) - \frac{1}{2L} \| \nabla f(x^k) \|_\infty^2 + \frac{\epsilon_{k}}{L} \sqrt{2L_1(f(x^0)-f(x^*))} - \frac{\epsilon_{k}^2}{2L} \\
& \leq f(x^k) - \frac{1}{2L} \| \nabla f(x^k) \|_\infty^2 + \epsilon_k\frac{\sqrt{2L_1}}{L}\sqrt{f(x^0)-f(x^*)}.
\end{align*}
Applying strong convexity (taken with respect to the $1$-norm), we get
$$
	f(x^{k+1}) - f(x^*) \le \bigg ( 1 - \frac{\mu_1}{L} \bigg ) \big [f(x^k) - f(x^*) \big ] + \epsilon_k\frac{\sqrt{2L_1}}{L}\sqrt{f(x^0)-f(x^*)},
$$
which implies
$$
\begin{aligned}
	f(x^{k+1}) - f(x^*) 
	& \le \bigg ( 1 - \frac{\mu_1}{L} \bigg )^k \big [f(x^0) - f(x^*) \big ] + \sum_{i = 1}^k \bigg ( 1 - \frac{\mu_1}{L} \bigg )^{k-i} \epsilon_i\frac{\sqrt{2L_1}}{L}\sqrt{f(x^0)-f(x^*)}  \\
	& = \bigg ( 1 - \frac{\mu_1}{L} \bigg )^k \bigg [f(x^0) - f(x^*) + \sqrt{f(x^0)-f(x^*)}A_k \bigg ],\\
\end{aligned}
$$
where 
$$
	A_k = \frac{\sqrt{2L_1}}{L}\sum_{i = 1}^k \bigg ( 1 - \frac{\mu_1}{L} \bigg )^{-i} \epsilon_i.
$$

\subsection{Additive error bound in terms of $L$}

By our additive error inequality, we have
\[
	| \nabla_{i_k} f(x^k) | + \epsilon_k \ge \| \nabla f(x^k) \|_\infty.
\]
Using this again in~\eqref{eq:additive} we get
\begin{align*}
	f(x^{k+1}) 
	& \le f(x^k) - \frac{1}{2L} \| \nabla f(x^k) \|_\infty^2 + \frac{\epsilon_{k}}{L}  \| \nabla f(x^k) \|_\infty - \frac{\epsilon_{k}^2}{2L} \\
	& \le f(x^k) - \frac{1}{2L} \| \nabla f(x^k) \|_\infty^2 + \frac{\epsilon_{k}}{L}  \big (  | \nabla_{i_k} f(x^k) | + \epsilon_k  \big ) - \frac{\epsilon_{k}^2}{2L} \\
	& = f(x^k) - \frac{1}{2L} \| \nabla f(x^k) \|_\infty^2 + \frac{\epsilon_k}{L} | \nabla_{i_k} f(x^k) | + \frac{\epsilon_k^2}{2L}.
\end{align*}
Further, from our basic progress bound that holds for any $i_k$ we have
\[
	f(x^*) \le f(x^{k+1}) \le f(x^k) - \frac{1}{2L} \bigg [ \nabla_{i_k} f(x^k ) \bigg ]^2 \le f(x^0) - \frac{1}{2L} \bigg [ \nabla_{i_k} f(x^k ) \bigg ]^2,
\]
which implies
\begin{align*}
 | \nabla_{i_k} f(x^k ) | & \le \sqrt{2L(f(x^0) - f(x^*))}.\\ 
\end{align*}
%If $i_k \not = i_{k-1}$, then $( x^k - x^* )^2_{i_k} = ( x^{k-1} - x^* )_{i_k}^2$. 
%Else, $i_k = i_{k-1}$ and
%%\[
%%	( x^k - x^* )^2_{i_k} \le ( x^{k-1} - x^* )_{i_k}^2.
%%\]
%\begin{align*}
%& ( x^k - x^* )^2_{i_k} \\
%	&~~= ( x^k - x^{k-1} + x^{k-1} - x^* )^2_{i_k} \\
%	&~~= (x^{k-1} - x^* )^2_{i_k} + 2( x^k - x^{k-1})_{i_k} (x^{k-1} - x^* )_{i_k} +  ( x^k - x^{k-1})_{i_k}^2\\
%	&~~= (x^{k-1} - x^* )^2_{i_k} - \frac{2}{L} \nabla_{i_k} f(x^{k-1}) (x^{k-1} - x^* )_{i_k} + \frac{1}{L^2} (\nabla_{i_k} f(x^{k-1}))^2\\
%	&~~= (x^{k-1} - x^* )^2_{i_k} - \frac{2}{L} (\nabla_{i_k} f(x^{k-1}) - \nabla_{i_k} f(x^*)) (x^{k-1} - x^* )_{i_k} + \frac{1}{L^2} (\nabla_{i_k} f(x^{k-1}))^2, \quad (\text{as } \nabla_{i_k} f(x^*) = 0)\\
%	&~~= (x^{k-1} - x^* )^2_{i_k} - \frac{2}{L} \big (\nabla_{i_k} f(x^* + h_{i_k}(x^{k-1}) e_{i_k}) - \nabla_{i_k} f(x^*) \big ) (x^* + h_{i_k} e_{i_k} - x^* )_{i_k} + \frac{1}{L^2} (\nabla_{i_k} f(x^{k-1}))^2\\
%	&~~\le (x^{k-1} - x^* )^2_{i_k} - \frac{2}{L^2} (\nabla_{i_k} f(x^{k-1}))^2 + \frac{1}{L^2} (\nabla_{i_k} f(x^{k-1}))^2, \quad (\text{Nesterov 2.1.8})\\
%	&~~= (x^{k-1} - x^* )^2_{i_k} - \frac{1}{L^2} (\nabla_{i_k} f(x^{k-1}))^2\\
%	&~~\le (x^{k-1} - x^* )^2_{i_k}.\\
%\end{align*}
%Applying this argument recursively gives us
%\[
%	( x^k - x^* )^2_{i_k} \le ( x^{0} - x^* )_{i_k}^2
%\]
and thus that
\begin{align*}
	f(x^{k+1}) &\le f(x^k) - \frac{1}{2L} \| \nabla f(x^k) \|_\infty^2 + \frac{\epsilon_k}{L} \sqrt{2L(f(x^0) - f(x^*))} + \frac{\epsilon_k^2}{2L}\\
	&= f(x^k) - \frac{1}{2L} \| \nabla f(x^k) \|_\infty^2 + \epsilon_k \sqrt{\frac{2}{L}}\sqrt{f(x^0) - f(x^*)} + \frac{\epsilon_k^2}{2L}.\\
\end{align*}
Applying strong convexity and applying the inequality recursively we obtain
%letting $\kappa = \big[ f(x^0) - f(x^*) \big]^{1/2}$ we get
%$$
%	f(x^{k+1}) - f(x^*) \le \bigg ( 1 - \frac{\mu_1}{L} \bigg ) \big [f(x^k) - f(x^*) \big ] + \sqrt{\frac{2}{L}} \epsilon_k \kappa + \frac{\epsilon_k^2}{2L},
%$$
%which implies
$$
\begin{aligned}
	f(x^{k+1}) - f(x^*) 
	& \le \bigg ( 1 - \frac{\mu_1}{L} \bigg )^k \big [f(x^0) - f(x^*) \big ] + \sum_{i = 1}^k \bigg ( 1 - \frac{\mu_1}{L} \bigg )^{k-i} \bigg ( \epsilon_i\sqrt{\frac{2}{L}} \sqrt{f(x^0) - f(x^*)} + \frac{\epsilon_i^2}{2L} \bigg ) \\
	& = \bigg ( 1 - \frac{\mu_1}{L} \bigg )^k \bigg [f(x^0) - f(x^*) + A_k \bigg ], \\
\end{aligned}
$$
where 
$$
	A_k = \sum_{i = 1}^k \bigg ( 1 - \frac{\mu_1}{L} \bigg )^{-i}\bigg ( \sqrt{\frac{2}{L}} \epsilon_i \sqrt{f(x^0) - f(x^*)} + \frac{\epsilon_i^2}{2L} \bigg ).
$$
Although uglier than the expression depending on $L_1$, this expression will tend to be smaller unless $\epsilon_k$ is not small.

%From~\citet{friedlander2011hybrid}: The bound on the error depends on both the conditioning of the problem (as determined by $L$ and $\mu$), and on the lower bound on the distance to the optimal function value $\pi_k$. We can allow a larger error in the gradient calculation the further the current iterate is from the optimal function value, but a more accurate calculation is needed to maintain the strong linear convergence rate as the iterates approach the solution. Similarly, if the problem is well conditioned so that the ratio $\mu/L$ is close to 1, a larger error in the gradient calculation is permitted, but for ill-conditioned problems, where $\mu/L$ is very small, we require a more accurate gradient calculation.

%%%%%%%%%%%%%%%%%%%%%%%%%%%%%%%%%%%%%%%%%%%%%%%%%%%

\section{Convergence Analysis of GS-$s$, GS-$r$, and GS-$q$ Rules}
\label{app:prox}

In this section, we consider problems of the form
\[
\min_{x \in \R^n} F(x) = f(x) + g(x) = f(x) + \sum_{i=1}^n g_i(x_i),
\]
where $f$ satisfies our usual assumptions, but the $g_i$ can be non-smooth. We first introduce some notation and state the convergence result,  and then show that it holds. We then show that the rate cannot hold in general for the GS-$s$ and GS-$r$ rules.

\subsection{Notation and basic inequality}

To analyze this case, an important inequality we will use is that the $L$-Lipschitz-continuity of $\nabla_i f$ implies that for all $x$, $i$, and $d$, we have
\begin{equation}
\label{eq:proxProgress}
\begin{split}
	F(x + d e_i) = f(x + d e_i) + g(x + d e_i) 
	&\le f(x) + \langle \nabla f(x), d e_i \rangle + \frac{L}{2} d^2 + g(x + d e_i) \\
	&= f(x) + g(x) + \langle \nabla f(x), d e_i \rangle + \frac{L}{2} d^2 + g_i(x_i + d) - g_i(x_i)\\
	&= F(x) + V_i(x,d), 
\end{split}
\end{equation}
where
\[
V_i(x,d) \equiv \langle \nabla f(x), d e_i \rangle + \frac{L}{2} d^2 + g_i(x_i + d) - g_i(x_i).
\]
Notice that the GS-$q$ rule is defined by
\[
i_k = \argmin{i}\{\min_d V_i(x,d)\}.
\]
We use the notation $d_i^k = \argmin{d} V_i(x^k,d)$ and we will use $d^k$ to denote the vector containing these values for all $i$. When using the GS-$q$ rule, the iteration is defined by
\begin{equation}
\label{eq:GSqUpdate}
\begin{aligned}
x^{k+1} & = x^k + d_{i_k}e_{i_k}\\
& = x^k + \argmin{d}\{V_{i_k}(x,d)\}e_{i_k}.
\end{aligned}
\end{equation}
In this notation the GS-$r$ rule is given by
\[
j_k = \argmax{i}|d_i^k|.
\]
Under this notation, we can show that coordinate descent with the GS-$q$ rule satisfies the bound
\begin{equation}
\label{eq:GSqRate}
F(x^{k+1}) - F(x^*) \leq  \left(1-\frac{\mu}{Ln}\right)[f(x^k)-f(x^*)].
 \end{equation}
% and $d^k$ is the (signed) difference between $x^k$ and the proximal-gradient step in all coordinates,
 %\[
 %d^k = x^k - \prox{\frac{1}{L}g} \left[ x_i^k - \frac{1}{L} \nabla f(x^k) \right],
 %\]
 %while $s^k$ is a particular element of the sub-differential $\partial g(x^k + d^k)$. 
 We show  this result by  showing that the GS-$q$ rule makes at least as much progress as randomized selection.
 % Finally, we also give a further bound that is sometimes faster than this bound, but relies on an assumption that is difficult to verify in practice.

 \subsection{GS-$q$ is at least as fast as random}
 
Our argument in this section follows a similar approach to~\citet{richtarik2014iteration}. In particular, combining~\eqref{eq:proxProgress} and~\eqref{eq:GSqUpdate} we have the following upper bound on the iteration progress
$$
\begin{aligned}
	F(x^{k+1}) 
	&\le F(x^k) + \min_{i\in\{1,2,\dots,n\}}  \left\{\min_{d \in \R}
		V_i(x^k,d)\right\}, \\
	&= F(x^k) + \min_{i\in\{1,2,\dots,n\}}  \left\{\min_{y \in \R^n}
		V_i(x^k,y - x^k)\right\}, \\
	&= F(x^k) + \min_{y \in \R^n} \left\{\min_{i\in\{1,2,\dots,n\}} 
		V_i(x^k,y - x^k)\right\}, \\
	&\le F(x^k) + \min_{y \in \R^n} \left\{\frac{1}{n} \sum_{i = 1}^n V_i(x^k,y - x^k)\right\} \\
	&= F(x^k) + \frac{1}{n} \min_{y \in \R^n} \left\{ \langle \nabla f(x^k), y-x^k \rangle + \frac{L}{2} \| y - x^k \|^2 +  g(y) - g(x^k) \right\} \\
%	&= f(x^k) + \bigg ( 1 - \frac{1}{n} \bigg ) g(x^k) + \langle \frac{1}{n}  \nabla f(x^k), y-x^k \rangle + \frac{L}{2n} \| y - x^k \|^2 +  \frac{1}{n} g(y).  \\
	&= \bigg ( 1 - \frac{1}{n} \bigg ) F(x^k) + \frac{1}{n} \min_{y \in \R^n} \left\{f(x^k) + \langle \nabla f(x^k), y-x^k \rangle + \frac{L}{2} \| y - x^k \|^2 +  g(y) \right\}.  \\
\end{aligned}
$$
From strong convexity of $f$, we have that $F$ is also $\mu$-strongly convex and that
\begin{align*}
f(x^k) & \leq f(y) - \langle \nabla f(x^k),y-x^k)\rangle - \frac{\mu}{2}\norm{y-x^k}^2,\\
F(\alpha x^* + (1-\alpha)x^k) & \leq \alpha F(x^*) + (1-\alpha)F(x^k) - \frac{\alpha(1-\alpha)\mu}{2}\norm{x^k-x^*}^2,
\end{align*}
for any $y \in \R^n$ and any $\alpha \in [0,1]$~\citep[see][Theorem~2.1.9]{Nes04b}.
Using these gives us
$$
\begin{aligned}
	& F(x^{k+1}) \\
	&~~~\le  \bigg ( 1 - \frac{1}{n} \bigg ) F(x^k) + \frac{1}{n} \min_{y \in \R^n} \left\{  f(y) - \frac{\mu}{2} \|y - x \|^2 + \frac{L}{2} \| y - x^k \|^2 +  g(y) \right\} \\
	&~~~=  \bigg ( 1 - \frac{1}{n} \bigg ) F(x^k) + \frac{1}{n} \min_{y \in \R^n} \left\{ F(y) + \frac{L - \mu}{2} \| y - x^k \|^2 \right\} \\
	&~~~\le  \bigg ( 1 - \frac{1}{n} \bigg ) F(x^k) + \frac{1}{n} \min_{\alpha \in [0,1] } \left\{F(\alpha x^* + (1-\alpha)x^k) + \frac{\alpha^2(L - \mu)}{2} \| x^k - x^* \|^2 \right\}\\
%	&\le  \bigg ( 1 - \frac{1}{n} \bigg ) F(x^k) + \frac{1}{n} \bigg [ \min_{\alpha \in [0,1] } \alpha F (x^*) + (1-\alpha)F(x^k)  - \frac{(\mu^f + \mu^g) \alpha (1-\alpha)}{2} \|x^k - x^* \|^2 + \frac{\alpha^2(L - \mu^f)}{2} \| x^k - x^* \|^2 \bigg ] \\
	&~~~\le  \bigg ( 1 - \frac{1}{n} \bigg ) F(x^k) + \frac{1}{n}  \min_{\alpha \in [0,1] }\left\{ \alpha F (x^*) + (1-\alpha)F(x^k)  + \frac{\alpha^2(L - \mu) \!-\! \alpha (1-\alpha)\mu}{2} \|x^k - x^* \|^2 \right\}\\
	&~~~\le  \bigg ( 1 - \frac{1}{n} \bigg ) F(x^k) + \frac{1}{n} \bigg [ \alpha^* F (x^*) + (1-\alpha^*)F(x^k) \bigg ] \quad\quad \left(\text{choosing } \alpha^* = \frac{\mu}{L}\in(0,1]\right) \\
	&~~~=  \bigg ( 1 - \frac{1}{n} \bigg ) F(x^k) + \frac{\alpha^*}{n} F (x^*) + \frac{(1-\alpha^*)}{n} F(x^k) \\
	&~~~=  F(x^k) - \frac{\alpha^*}{n} [F(x^k) - F(x^*)].
\end{aligned}
$$
Subtracting $F(x^*)$ from both sides of this inequality gives us
$$
	F(x^{k+1}) - F(x^*) \le \bigg (1 - \frac{\mu}{nL } \bigg )[F(x^k) - F(x^*)].
$$

\subsection{Lack of progress of the GS-$s$ rule}

We now show that the rate $(1-\mu/Ln)$ cannot hold  for the GS-$s$ rule. We do this by constructing a problem where an iteration of the GS-$s$ method does not make sufficient progress.
In particular, consider the bound-constrained problem
\[
	\min_{x \in C} f(x) = \frac{1}{2} \| Ax - b \|^2_2,
\]
where $C = \{ x ~:~ x \ge 0 \}$, and
\[
	A = \begin{pmatrix}   1 & 0 \\
					0 & 0.7  
	\end{pmatrix},
	\quad 
		b = \begin{pmatrix} -1 \\ -3 \end{pmatrix},
	\quad
	x^0 = \begin{pmatrix} 1 \\ 0.1 \end{pmatrix},
	\quad
	x^* = \begin{pmatrix} 0 \\ 0 \end{pmatrix}.
\]
We thus have that
\begin{align*}
f(x^0) & = \frac{1}{2}((1+1)^2 + (.07+3)^2) \approx 6.7\\
f(x^*) & = \frac{1}{2}((-1)^2 + (-3)^2) = 5\\
\nabla f(x^0) & = A^T(Ax_0 - b) \approx \begin{pmatrix} 2.0 \\ 2.1\end{pmatrix}\\
\nabla^2 f(x) & = A^TA = \begin{pmatrix}   1 & 0 \\
					0 & 0.49  
	\end{pmatrix}.
	\end{align*}
	The parameter values for this problem are
	\begin{align*}
  n &= 2 \\                     
  \mu &= \lambda_{min} = 0.49 \\
  L &= \lambda_{max} = 1 \\
  \mu_1 &= \left(\frac{1}{\lambda_1} + \frac{1}{\lambda_2}\right)^{-1} = 1 + \frac{1}{0.49} \approx 0.33,
\end{align*}
where the $\lambda_i$ are the eigenvalues of $A^T A$, and $\mu$ and $\mu_1$ are the corresponding strong-convexity constants for the $2$-norm and $1$-norm, respectively.

The proximal operator of the indicator function is the projection onto the set $C$, which involves  setting negative elements to zero. Thus, our iteration update is given by 
\[
	x^{k+1} = \prox_{\delta_C} [x^k - \frac{1}{L}\nabla_{i_k} f(x^k)e_{i_k}] = \max(x^k - \frac{1}{L}\nabla_{i_k} f(x^k)e_{i_k}, 0),
\]
For this problem, the GS-$s$ rule is given by
\[
	i = \argmax{i} | \eta_i^k |,
\]
where 
\[
	\eta_i^k = \begin{cases}
	\nabla_i f(x^k), & \text{if } x_i^k \not = 0 \text{ or } \nabla_i f(x^k) < 0 \\
	0  , & \text{otherwise}
	\end{cases}.
\]
Based on the value of $\nabla f(x^0)$, the GS-$s$ rule thus chooses to update coordinate 2, setting it to zero and obtaining
\[
f(x^1) = \frac{1}{2}((1+1)^2 + (-3)^2) = 6.5.
\]
Thus we have
\[
\frac{f(x^1)-f(x^*)}{f(x^0)-f(x^*)} \approx \frac{6.5-5}{6.7-5} \approx 0.88,
\]
which we contrast with the bounds of
\begin{align*}
\left(1-\frac{\mu}{Ln}\right) & = \left(1 - \frac{0.49}{2}\right) \approx 0.76,\\
 \left(1 - \frac{\mu_1}{L}\right) & \approx (1-0.33) = 0.67.
\end{align*}
Thus, the GS-$s$ rule does not satisfy either bound. On the other hand, the GS-$r$ and GS-$q$ rules are given in this context by
\[
	i_k = \argmax{i}\left|\max\left(x^k - \frac{1}{L}\nabla_i f(x^k)e_i, 0\right) - x^k\right|,
\]
and thus both these rules choose to update coordinate 1, setting it to zero to obtain $f(x^1) \approx 5.2$ and a progress ratio of
\[
\frac{f(x^1)-f(x^*)}{f(x^0)-f(x^*)} \approx \frac{5.2-5}{6.7-5} \approx 0.12,
\]
which clearly satisfies both bounds.

\subsection{Lack of progress of the GS-$r$ rule}

We now turn to showing that the GS-$r$ rule does not satisfy these bounds in general. It will not be possible to show this for a simple bound-constrained problem since the GS-$r$ and GS-$q$ rules are equivalent for these problems. Thus, we consider the following $\ell_1$-regularized problem
\[
	\min_{x \in \R^2} \frac{1}{2} \| Ax - b \|^2_2 + \lambda \| x \|_1 \equiv F(x).
\]
We use the same $A$ as the previous section, so that $n$, $\mu$, $L$, and $\mu_1$ are the same. However, we now take
\[
		b = \begin{pmatrix} 2 \\ -1 \end{pmatrix},
		\quad 
		x_0 = \begin{pmatrix} 0.4 \\ 0.5 \end{pmatrix},
	\quad
	x_* = \begin{pmatrix} 1 \\ 0 \end{pmatrix},\quad 
	\lambda = 1,
\]
so we have
\begin{align*}
f(x_0)   \approx 3.1, \quad
f(x_*)  = 2
\end{align*}
The proximal operator of the absolute value function is given by the soft-threshold function, and our coordinate update of variable $i_k$ is given by
\[
	x_{i_k}^{k+1} = \prox_{\lambda|\cdot|} [x_{i_k}^{k+\frac12}] = \sgn{x_{i_k}^{k+\frac12}}\cdot \max(x_{i_k}^{k+\frac12} - \lambda/L, 0),
\]
where we have used the notation
\[
	x_i^{k+\frac12} = x_i^k - \frac{1}{L} \nabla_i f(x^k)e_i.
\]
The GS-$r$ rule is defined by
\[
i_k = \argmax{i}|d_i^k|,
\]
where $d_i^k = \prox_{\lambda |\cdot|} [x_i^{k+\frac12}] - x_i^k$ and in this case
\[
	d^0 = \begin{pmatrix} 0.6 \\ -0.5 \end{pmatrix}.
\]
Thus, the GS-$r$ rule chooses to update coordinate $1$. After this update the function value is
\[
	F(x^1) \approx 2.9,
\]
so the progress ratio is
\[
\frac{F(x^1)-F(x^*)}{F(x^0)-F(x^*)} \approx \frac{2.9 - 2}{3.1 - 2} \approx 0.84.
\]
However, the bounds suggest faster progress ratios of
\[
	\bigg (1 - \frac{\mu}{Ln} \bigg ) \approx 0.76,
\]
\[
	\bigg (1 - \frac{\mu_1}{L} \bigg ) \approx 0.67,
\]
so the GS-$r$ rule does not satisfy either bound. In contrast, in this setting the GS-$q$ rule chooses to update coordinate $2$ and obtains $F(x^1) \approx 2.2$, obtaining a progress ratio of
\[
\frac{F(x^1)-F(x^*)}{F(x^0)-F(x^*)} \approx \frac{2.2 - 2}{3.1 - 2} \approx 0.16,
\]
which satisfies both bounds by a substantial margin. Indeed, we used a genetic algorithm to search for a setting of the parameters of this problem (values of $x^0$, $\lambda$, $b$, and the diagonals of $A$) that would make the GS-$q$ rule not satisfy the bound depending on $\mu_1$, and it easily found counter-examples for the GS-$s$ and GS-$r$ rules but was not able to produce a counter example for the GS-$q$ rule.

\section{Runtime Experiments}
\label{app:RT}

\begin{figure}
\centering
\includegraphics[width=.5\textwidth]{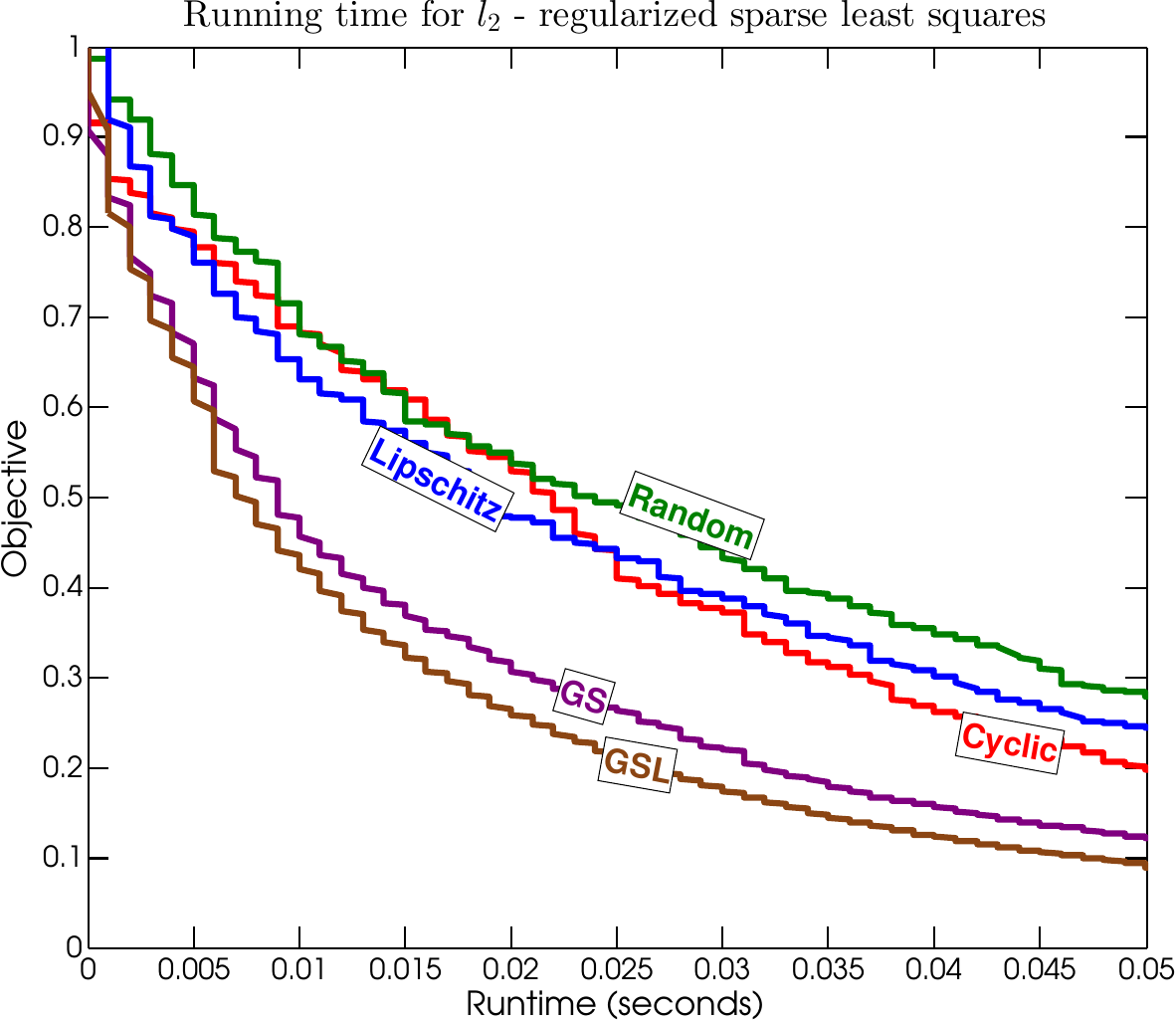}
\caption{Comparison of coordinate selection rules for $\ell_2$-regularized sparse least squares.}
\label{fig:time}
\end{figure}

In Figure~\ref{fig:time} we plot the objective against the runtime for the $\ell_2$-regularized sparse least squares problem from our experiments. Although runtimes are very sensitive to exact implementation details and we believe that more clever implementations than our naive Python script are possible, this figure does show that the GS and GSL rules offer benefits in terms of runtime with our implementation and test hardware.

\end{appendices}

\bibliography{bib}
\bibliographystyle{abbrvnat}

\end{document}